\documentclass[12pt]{article}
\usepackage{latexsym,amssymb,amsmath}

\newtheorem{thm}{Theorem}

\newtheorem{lemma}[thm]{Lemma}

\newcommand{\qed}{\hfill$\Box$}

\def\eps{\varepsilon}

\def\l{\lambda}
\def\lam{\lambda}
\def\Lam{\Lambda}

\def\a{\alpha}

\def\ed#1{{\mathbf 1}_{\{#1\}}}

\def\E{{\mathbb E}}
\def\ee{\textrm{\tiny E}}
\newcommand{\reals} {{\mathbb R}}      
\def\ww{\textrm{\tiny W}}
\def\OO{{\mathcal O}}

\def\xx{{\rm x}}

\def\N{{{\rm I}\!{\bf N}}}

\def\il{\left<}
\def\ir{\right>}

\def\e0{e(Q_{0,d})}

\title{Tractability of Multi-Parametric Euler and Wiener Integrated Processes}

\author{ M. A. Lifshits, A. Papageorgiou, H. Wo\'zniakowski }

\begin{document}
\maketitle
\date
\vskip 4pc

\begin{abstract}

\vskip 1pc

We study average case approximation of Euler and Wiener integrated processes
of~$d$ variables which are almost surely $r_k$-times continuously 
differentiable with 
respect to the $k$-th variable and $0\le r_k\le r_{k+1}$. Let $n(\eps,d)$ 
denote the minimal number of continuous 
linear functionals which is needed to 
find an algorithm that uses $n$ such functionals and whose average case error
improves the average case error of the zero algorithm by a factor $\eps$. 
Strong polynomial tractability means that there are nonnegative numbers $C$ and 
$p$ such that 
$$
n(\eps,d)\le C\eps^{-p}\ \ \ \ \ 
\mbox{for all}\ \  \ d\in\N=\{1,2,\dots\},\ \ \mbox{and}\ \ \eps\in(0,1).
$$
We prove that the Wiener process is much more difficult to approximate than the 
Euler process. Namely, strong polynomial tractability holds for the Euler case iff 
$$
   \liminf_{k\to \infty}\ \frac {r_k}{\ln\,k}> \frac1{2\ln\,3},
$$
whereas it holds for the Wiener case iff
$$
   \liminf_{k\to \infty}\ \frac {r_k}{k^s}> 0\ \ \ \ \ \mbox{for some}\ \
  s>\tfrac12.
$$
Other types of tractability are also studied. 
\end{abstract}

\section{Introduction}\label{s1}

Tractability of multivariate problems has been recently an active research
area. The reader may see the current progress on tractability 
in~\cite{NW08,NW10,NW12}. Tractability has been studied in various
settings and for various error criteria.

This field deals with problems defined on spaces of $d$-variate 
functions. For many practical computational problems $d$ is large. This
holds for problems in mathematical finance, statistics and physics.
We usually want to solve multivariate problems to within an error
threshold $\eps$ by algorithms that use finitely many function values
or, more generally, finitely many continuous linear functionals. 
Let $n(\eps,d)$ be the information complexity or shortly the complexity,
denoting the minimal number of function values or continuous linear functionals 
that are needed to find an algorithm approximating the solution of 
a multivariate problem to within $\eps$.  

Many multivariate problems suffer from the \emph{curse of dimensionality}.
That is, $n(\eps,d)$ is exponentially large in $d$. 
One of the goals of tractability is to determine under which conditions
the curse of dimensionality is \emph{not} present. Even more, we would like
to have the complexity bounded by some non-exponential function of
$d$ and $\eps^{-1}$. In particular, we have 
\begin{itemize}
\item \emph{weak} tractability if the complexity is \emph{not}
exponential in $d$ or $\eps^{-1}$, 
\item \emph{quasi-polynomial}
tractability if the complexity is of order 
$\exp(\,t\,(1+\ln\,d)(1+\ln\,\eps^{-1}))$,
\item \emph{polynomial} tractability if the complexity is of order
$d^{\,q}\,\eps^{-p}$,
\item \emph{strong} polynomial tractability if the complexity
is of order $\eps^{-p}$.
\end{itemize}

All bounds above hold for all $d$ and all $\eps\in(0,1)$ 
with the parameters $t,q,p$ and the factors multiplying the
corresponding complexity bounds independent of $d$ and 
$\eps^{-1}$.

The strong polynomial tractability is the most challenging property. 
If this holds then the complexity has a bound independent of $d$.
One may think that this property may hold only for trivial problems. 
Luckily, the opposite is true.

The curse of dimensionality often holds for multivariate problems for
which all variables and groups of variables play the same role. 
One way to vanquish the curse is to shrink the class of functions by 
introducing weights that monitor the influence of successive
variables and groups of variables. For sufficiently fast decaying
weights we not only vanquish the curse but even obtain strong
polynomial tractability; a survey of such results may be found 
in~\cite{NW08,NW10,NW12}.

The other way to vanquish the curse is by increasing the 
smoothness of functions with respect to successive variables. 
This approach was taken recently in~\cite{PW10}. It was done for 
multivariate approximation defined over Korobov spaces in the worst
case setting. The current paper can be viewed as a continuation 
of~\cite{PW10}. We consider multivariate approximation but now 
in the average case setting with the normalized error criterion.
This error criterion is defined as follows. 
We first take the zero algorithm and find its average case error
for multivariate approximation; this is called the initial error.
The initial error tells us how the problem scales
and what can be achieved without sampling the functions.
The normalized error criterion means that we want to improve
the initial error by a factor~$\eps$. We analyze algorithms that
use arbitrary continuous linear functionals. We stress 
that the same results hold for algorithms that use 
only function values. This is due to general relations
between these two classes of algorithms established in~\cite{HWW08}
and in Chapter~24 of~\cite{NW12}. 

In this paper we analyze two multivariate approximation problems 
defined for the Euler and Wiener integrated processes, whereas
in~\cite{LPW2} we consider average case approximation for 
general non-homogeneous tensor products.  
More precisely, here we take the space of continuous real functions
defined on the $d$-dimensional unit cube $[0,1]^d$.
We stress that $d$ can be an arbitrary positive integer, however, our 
emphasis is on large $d$.  We equip this space with a zero-mean Gaussian 
measure whose covariance kernel is denoted by $K_d$. We study two such
kernels. The first one is $K_d=K_d^{\ee}$ for the Euler integrated process, 
whereas the second one $K_d=K_d^{\ww}$ is for the Wiener integrated
process. These processes are precisely defined  in the next section. 
Here we only mention that for both of them we know that almost surely the
functions are $r_k$ times continuously differentiable with respect to 
the $k$-th variable for $k=1,2,\dots,d$. 

The information complexity is then denoted
by $n^{\ee}(\eps,d)$ and $n^{\ww}(\eps,d)$ for the Euler and Wiener 
integrated processes, respectively. Obviously, it depends on the the smoothness 
parameters~$\{r_k\}$. Our main goal in this paper is to find necessary
and sufficient conditions in terms of~$\{r_k\}$ such that the four
notions of tractability are satisfied. 

We now briefly describe the results obtained in this paper.  
For both processes we prove that weak tractability holds iff
$\lim_{k\to\infty}r_k=\infty$. 
Otherwise, if $r=\lim_{k\to\infty}r_k<\infty$ then 
we have the curse of dimensionality. This means that if all
$r_k\le r<\infty$ then both $n^{\ee}(\eps,d)$ and $n^{\ww}(\eps,d)$ 
depend exponentially on $d$ and this holds for all
$\eps\in(0,1)$. Hence, the function $n^{\xx}(\cdot,d)$ is discontinuous 
at $1$. Indeed, $n^{\xx}(1,d)=0$ although for $\eps$ pathologically close 
to one $n^{\rm x}(\eps,d)$ depends exponentially on $d$; here 
$\xx \in\{ \textrm{\sc E,\sc W} \}$.  

We stress that weak tractability does not depend on the rate of
convergence of $r_k$ to infinity. However, if we want to obtain 
other types of tractability we must require a certain
convergence rate for the $r_k$, although the rate is different 
for the Euler and the Wiener case. 
For simplicity, let us consider 
$$
   r_k=\left\lceil 1+ a\,\ln(1+k\ln\,k)\right\rceil
$$
for some positive number $a$. Then for the Euler case we have:
\begin{itemize}
\item $a<\frac1{2\,\ln\,3}$\ \ \  no quasi-polynomial tractability,
\item $a=\frac1{2\,\ln\,3}$\ \ \  quasi-polynomial tractability but no
  polynomial tractability,          
\item $a>\frac1{2\,\ln\,3}$\ \ \  strong polynomial tractability.
\end{itemize}

For the Wiener case we have to assume much more since for $r_k$  given
above only weak tractability holds. For
$$
  r_k=\left\lceil k^s\,\ln^2(1+k)\right\rceil
$$
we have
\begin{itemize}
\item $s<\frac1{2}$\ \ \  no quasi-polynomial tractability,
\item $s=\frac1{2}$\ \ \  quasi-polynomial tractability but no
  polynomial tractability,          
\item $s>\frac1{2}$\ \ \  strong polynomial tractability.
\end{itemize}

For general $\{r_k\}$, we prove that quasi-polynomial tractability holds
iff
\begin{eqnarray*}
\mbox{For the Euler case}\ &:&\ \ \ \limsup_{d\to\infty}\ \frac1{\ln\,d}
\,\sum_{k=1}^d  (1+ r_k)\,3^{-2r_k}<\infty,\\
\mbox{For the Wiener case:}\ &:&\ \ \  \limsup_{d\to\infty}\ \frac1{\ln\,d}
\,\sum_{k=1}^d (1+r_k)^{-2}\,\max(1,\ln\,r_k)<\infty.
\end{eqnarray*}
Furthermore, for both processes polynomial tractability is equivalent
to strong polynomial tractability and holds iff
\begin{eqnarray*}
\mbox{For the Euler case}\ &:&\ \ \ 
a_{\ee}:=\liminf_{d\to\infty}\ \frac{r_k}{\ln\,k}>\frac1{2\,\ln\,3},\\
\mbox{For the Wiener case:}\ &:&\ \ \ 
\liminf_{d\to\infty}\ \frac{r_k}{k^{s}}>0
\ \ \ \ \ \mbox{for some}\ \ \ s>\tfrac12.
\end{eqnarray*}
We also study 
the exponent $p^{\,\rm str-avg-x}$ of strong polynomial tractability 
which is defined as the infimum of $p$ for which the complexity is of order
$\eps^{-p}$. For the Euler case we have 
$$
p^{\,\rm str-avg-\ee}
  =\max\left(\frac2{2r_1+1},\frac2{2 a_{\ee}\,\ln\,3 -1}\right)
$$
For the Wiener case and $r_k=k^s$ for some $s>\tfrac12$ we have 
$$
   \max\left(\frac 2{2r_1+1},\frac2{2s-1}\right)\le
   p^{\,\rm str-avg-\ww}\le \max\left(\frac2{2s-1},3\right).
$$
Hence, for $s\in(\tfrac12,\tfrac56]$ we know that
$$
   p^{\,\rm str-avg-\ww}=\frac2{2s-1}\, ;
$$
otherwise our bounds are too weak to provide the exact value of 
the exponent.

Our results solve a special case of Open Problem~11 in \cite{NW08}, where
$r_{d,k}=r_k$, $k=1,\dots,d$, and, with slightly modified proofs, 
they also solve Open Problem~10 in \cite{NW08}.

The Euler, Wiener and other univariate integrated processes can be 
characterized as follows.  Consider
\[
   X^{\a}(t):=
   (-1)^{\a_1+\dots+\a_r} \underbrace{
\int_{\a_r}^t \int_{\a_{r-1}}^{t_{r-1}} \dots \int_{\a_1}^{t_1} W(s)\ 
{\rm d}s {\rm d}t_1  \dots {\rm d}t_{r-1}
}_{ r\ {\rm times}},
\qquad  \qquad 0\le t\le 1,
\]
where $\a=(\a_1,\a_2,\dots, \a_r)$ is a multi-index with components 
$\a_i\in\{ 0,1\}$ for $i=1,2,\dots,r$,
$W(s)$ is the standard Wiener process for $0\le s\le 1$, 
and $r\in\N$. Then
$X^{(0,\dots,0)}$ is the 
integrated Wiener process and $X^{(1,0,1,0,1,\dots)}$ is the 
integrated Euler process. It is an open problem to 
consider the integrated processes
resulting from the different values of the multi-index $\a$ 
and to compare the necessary and sufficient conditions
on $\{ r_k\}$ for weak, quasi-polynomial and polynomial tractability, 
respectively, with those obtained for the Euler and Wiener processes.
In particular, it seems of interest to verify whether the Euler
process is the easiest and the Wiener process is the most difficult
among all of these $2^r$ processes. 

The paper is organized as follows. In Section 2 we present the precise
definitions of the average case approximation problem, the Euler and 
Wiener integrated processes and tractability notions. 
In Section 3 we present results for the Euler and  in Section 4 for the 
Wiener integrated processes. The proofs of three theorems are presented 
in Sections 5 to 7.

\section{Preliminaries}\label{s2}

In this section we precisely define the Euler and Wiener processes,
multivariate approximation in the average case setting, and we cite 
known results that will be needed for our analysis.

\subsection{Euler and Wiener Processes}\label{s2.1}

Let $F_d=C([0,1]^d)$ be the space of real continuous functions defined 
on $[0,1]^d$ with the max norm,
$$
\|f\|_{F_d}=\max_{x\in [0,1]^d}|f(x)|\ \ \ \ \ 
\mbox{for all}\ \ \ f\in F_d.
$$
We equip the space $F_d$ with a zero-mean Gaussian measure $\mu_d$
defined on Borel sets of $F_d$. The covariance kernel $K_d$ related to
 $\mu_d$ is defined by
$$
  K_d(x,t)=\int_{F_d}f(x)\,f(t)\,\mu_d({\rm d}f)\ \ \ \ \  
  \mbox{for all}\ \ \ x,t\in [0,1]^d.
$$
We refer to \cite{L96} for extensive theory of Gaussian measures in linear 
spaces and their covariance kernels. 

By $\{r_k\}$ we mean a sequence of non-negative non-decreasing integers 
$$
   0\le r_1\le r_2\le\cdots\le r_d\le \cdots\ .
$$ 
The Euler and Wiener integrated processes differ in the choice 
of the covariance kernel $K_d$.
Our presentation of the Euler integrated processes is based
on~\cite{CH01} and~\cite{GHT03}.
The Wiener integrated process is more standard and can be found in
many books and papers.
\begin{itemize}
\item
\emph{Euler integrated process.}\ \ 
We now have $K_d=K_d^{\ee}$ given by  
$$
K_d^{\ee}(x,y)=\prod_{k=1}^d K_{1,r_k}^{\ee}(x_k,y_k)\ \ \ 
\mbox{for all}\ \ \ x,y\in[0,1]^d,
$$
where
$$
K_{1,r}^{\ee}(x,y)=
\int_{[0,1]^r}\min(x,s_1)\,\min(s_1,s_2)\,\cdots\,\min(s_r,y)\,{\rm d}s_1
\,{\rm d}s_2\cdots{\rm d}s_r
$$
for all $x,y\in [0,1]$. This kernel is equal to
$$
K_{1,r}^{\ee}(x,y)=
(-1)^{r+1}\frac{2^{2r}}{(2r+1)!}
\bigg(E_{2r+1}(\tfrac12|x-y|)-E_{2r+1}(\tfrac12(x+y))\bigg)
$$
for all $x,y\in [0,1]$. Here, $E_n$ is the $n$-th degree Euler 
polynomial which can be defined as the coefficient of the generating 
function
$$
\frac{2\exp(x\,t)}{\exp(x)+1}=\sum_{n=0}^{\infty} E_n(x)\frac{t^n}{n!}\
\ \ \ \ \mbox{for all}\ \ \ x,t\in\reals.
$$
In particular, we have $E_0=1$, $E_1(x)=x-\tfrac12$ and
$E_2(x)=x^2-x$.

The process is called \emph{Euler} due to the fact that the 
covariance  kernel is expressed by Euler polynomials.  
\item
\emph{Wiener integrated process.}\ \ We now have $K_d=K_d^{\ww}$ given by 
$$
  K_d^{\ww}(x,y)=\prod_{k=1}^d K_{1,r_k}^{\ww}(x_k,y_k)\ \ \ 
  \mbox{for all}\ \ \ x,y\in[0,1]^d,
$$
where
$$
  K_{1,r}^{\ww}(x,y)=\int_0^{\min(x,y)}
  \frac{(x-u)^r}{r!}\,
  \frac{(y-u)^r}{r!}\,{\rm d}u
  =\int_0^1\frac{(x-u)_+^r}{r!}\,
  \frac{(y-u)_+^r}{r!}\,{\rm d}u
$$
for all $x,y\in [0,1]$ and with the standard notation $t_+=\max(t,0)$. 
\end{itemize}

Let us stress that the univariate Euler and Wiener processes are close
relatives since they emerge from very similar 
integration schemes. Indeed, 
let $W(t), t\in[0,1]$, be a standard Wiener process, i.e. a Gaussian random
process with zero mean and covariance
\[ 
K_{1,0}^\ee(s,t)=K_{1,0}^\ww(s,t)=\min(s,t).
\]
Consider two sequences of integrated processes $X_r, Y_r$ on $[0,1]$ 
defined by $X_0=Y_0=W$, and for $r=0,1,2,\dots$ 
\begin{eqnarray*}
  X_{r+1}(t)&=& \int_0^t X_{r}(s) {\rm d} s  
  \\
  Y_{r+1}(t)&=& \int_{1-t}^1 Y_{r}(s) {\rm d} s, 
\end{eqnarray*}
Then the covariance kernel of $X_r$ is $K_{1,r}^\ww$ 
while the covariance kernel of
$Y_r$ is $K_{1,r}^\ee$. Clearly, $X_r$ and $Y_r$ have the same smoothness
properties. That is why 
different tractability results are surprising.

On the other hand, there are some differences between the two processes.
The Gaussian measure $\mu_d$ on $F_d$ corresponding 
to the covariance kernel $K_d^{\ee}$ or 
$K_d^{\ww}$ is concentrated on functions that are almost surely $r_k$-times
continuously differentiable 
with respect to the $k$-th variable for $k=1,2,\dots,d$, and satisfy certain 
boundary conditions which are different for the Euler and Wiener cases. 

For the Euler case, we have
\begin{equation}\label{E1}
\frac{\partial^{k_1+k_2+\cdots +k_d}}{\partial\,x_1^{k_1}\,\partial\,
x_2^{k_2}\cdots\partial\,x_d^{k_d}}\, f(x) = 0
\end{equation}
if for some $i$ we have $x_i=0$ and $k_i$ is even, or
$x_i=1$ and $k_i$ is odd. Here, $k_i=0,1,\dots,r_i$.

For the Wiener case, we have
\begin{equation}\label{W1}
\frac{\partial^{k_1+k_2+\cdots +k_d}}{\partial\,x_1^{k_1}\,\partial\,
x_2^{k_2}\cdots\partial\,x_d^{k_d}}\, f(x) = 0
\end{equation}
if one of the components of $x$ is zero. As before,
$k_i=0,1,\dots,r_i$.

To see the difference between~\eqref{E1} and~\eqref{W1} more
explicitly, we take $d=1$. Then for the Euler case 
for all  $0\le k\le r_1$ we have
$$
f^{(k)}(0)=0 \ \ \ \mbox{if $k$ is even \ \ \ and}\ \ \
\ \ f^{(k)}(1)=0 \ \ \mbox{if $k$ is odd},
$$
whereas for the Wiener case we have
$$
f^{(k)}(0)=0\ \ \ \mbox{for all $k\le r_1$}.
$$

Finally, note that Nazarov and Nikitin studied in \cite{Naz,NazNik} a slightly
different version $W_r^N$ of Euler integrated process. 
The processes $W_r^N$ and 
$W_r^\ee$ coincide for even $r$ but  $W_r^N(t)=W_r^\ee(1-t)$ for odd $r$. 
The covariance spectra of both processes are the same but the boundary 
conditions are 
different. Since the spectra are the same, the tractability results 
for the Nazarov and Nikitin process are the same as for the Euler
process.

\subsection{Multivariate Approximation}\label{s2.2}

Multivariate approximation is defined by 
the embedding ${\rm APP}_d:F_d\to L_2$ given by 
$${\rm APP}_df=f\ \ \ \ \ \mbox{for all}\ \ \ f\in F_d.
$$
Here, $L_2=L_2([0,1]^d)$ is the standard $L_2$ space with
the norm 
$$
\|f\|_{L_2}=\bigg(\int_{[0,1]^d}f^{\,2}(t)\,{\rm d}t\bigg)^{1/2}.
$$ 
We approximate functions $f$ from $F_d$ by algorithms $A_n$ 
that use $n$ function values or arbitrary 
continuous linear functionals. We only consider the case of arbitrary
continuous functionals since it is known that the results are the same 
for function values, see~\cite{HWW08} and Chapter~24 of~\cite{NW12}.  
In the average case setting,
without essential loss of
generality, see e.g.,~\cite{TWW88} as well as~\cite{NW08}, 
we can restrict ourselves to
linear algorithms $A_n$ of the form
$$
A_n(f)=\sum_{j=1}^nL_j(f)\,g_j \ \ \ \ \  
\mbox{with} \ \ \ \ \ L_j\in F^{\,*}_d,\ \ g_j\in L_2.
$$
The average case error of $A_n$ is defined as
$$
e^{\rm avg}(A_n)=\bigg(\int_{F_d}\|{\rm APP}_df-A_n(f)\|_{L_2}^2\,
\mu_d({\rm   d}f)\bigg)^{1/2}, 
$$
where $\mu_d$ is a zero-mean Gaussian measure with a covariance kernel
$K_d$ and  
$$
\|{\rm APP}_df-A_n(f)\|_{L_2}^2=\int_{[0,1]^d}(f(t)-A_n(f)(t))^2\,{\rm d}t.
$$
Then $\nu_d=\mu_d{\rm APP}_d^{-1}$ is a zero-mean Gaussian measure 
defined on Borel sets of $L_2$. The covariance operator 
$C_{\nu_d}:L_2\to L_2$  of $\nu_d$ is given by
$$
  C_{\nu_d}f= \int_{[0,1]^d}K_d(\cdot,t)\,f(t)\,{\rm d}t\ \ \ \ \ 
  \mbox{for all}\ \ \ f\in L_2.
$$
The operator $C_{\nu_d}$ is a self-adjoint, nonnegative definite,
and has finite trace. Let 
$(\l_{d,j},\eta_{d,j})_{j=1,2,\dots}$ denote its eigenpairs 
$$
  C_{\nu_d}\eta_{d,j}=\lambda_{d,j}\,\eta_{d,j} \ \ \ \ \
  \mbox{with} \ \ \ \lambda_{d,1}\ge\lambda_{d,2}\ge\cdots,
$$
and 
$$
  \sum_{j=1}^\infty\lambda_{d,j}=\int_{[0,1]^d}K_d(t,t)\,{\rm d}t<\infty.
$$

For a given $n$, it is well known that the algorithm $A_n$ that
minimizes the average case error is of the form
\begin{equation}\label{optalg}
   A_n(f)=\sum_{j=1}^n\il f,\eta_{d,j}\ir_{L_2}\eta_{d,j},
\end{equation}
and its average case error is
\begin{equation}\label{avgerror}
    e^{\rm avg}(A_n)=\bigg(\sum_{j=n+1}^\infty\lambda_{d,j}\bigg)^{1/2}.
\end{equation}
For $n=0$ we obtain the zero algorithm $A_0=0$. Its average case
error is called the initial error, and is given by the square-root of
the trace of the operator $C_{\nu_d}$, i.e., by~\eqref{avgerror} with
$n=0$.

We now define the average case information complexity 
$n(\eps,d)$ as the minimal~$n$ for which there is an algorithm whose 
average case error reduces the initial error by a factor~$\eps$,
\begin{equation}\label{infocomp}
   n(\eps,d)
   =\min\bigg\{\,n\ \bigg| \ \ \sum_{j=n+1}^\infty\lambda_{d,j}\le
   \eps^2\,\sum_{j=1}^\infty\lambda_{d,j}\bigg\}.
\end{equation}

{}From~\eqref{infocomp} it is clear that all notions of tractability
depend only on the eigenvalues $\lambda_{d,j}$. Therefore the more we
know about the eigenvalues $\lambda_{d,j}$ the more we can say about
various notions of tractability.

\subsection{Eigenvalues for the Euler and Wiener Integrated 
  Processes}\label{s3.1}

For both processes the corresponding covariance kernel is of product
form. Therefore the eigenvalues for the $d$-variate case are products
of the eigenvalues of the univariate cases which depend on the
smoothness parameters $r_k$ for $k=1,2,\dots,d$. That is, if we denote
   by $\lambda_{d,j}^\xx$'s  
the eigenvalues of the Euler integrated process,
$\xx={\rm E}$, or  
the eigenvalues of the Wiener integrated
process, $\xx={\rm W}$, then
$$
   \{\lambda_{d,j}^\xx\}_{j=1,2,\dots}
   = \left\{\lambda_{j_1,r_1}^\xx \lambda_{j_2,r_2}^\xx \dots
   \lambda_{j_d,r_d}^\xx \right\}_{j_1,j_2,\dots,j_d=1,2\dots},
$$ 
with the $\lambda_{j_k,r_k}^\xx$'s denoting the eigenvalues of the
univariate case with smoothness $r_k$.

For the Euler case, the $\lambda_{j_k,r_k}^{\ee}$'s are the eigenvalues of
the operator
$$
(C_{1,r_k}^{\,\ee}f)(x)=\int_0^1K^{\ee}_{1,r_k}(x,t)\,f(t)\,\,{\rm d}t.
$$
By successive differentiation of this equation with respect to $x$
and using the properties of the kernel $K_{1,r_k}^{\ee}$, it is easy
to show that the eigenvalues of $C_{1,r_k}^{\,\ee}$ satisfy
the Sturm-Liouville problem
\begin{equation}\label{SL}
   \lambda\,f^{(2r_k+2)}(x)=(-1)^{r_k+1}f(x)\ \ \ \ \ \mbox{for all}
   \ \ \ x\in(0,1),
\end{equation}
with the boundary conditions
$$
   f(t_0)
   =f^{\prime}(t_1)=f^{\prime\prime}(t_2)=\cdots=f^{(2r_k+1)}(t_{2r_k+1})
   =0,
$$
where $t_i=0$ for even $i$ and $t_i=1$ for odd $i$. For the Euler case, 
we know the eigenvalues exactly, see \cite{CH01} and \cite{GHT03},
and they are equal to
\begin{equation}\label{Eulereig}
   \lambda_{j,r_k}^{\ee}=\left(\frac{1}{\pi(j-1/2)}\right)^{2r_k+2}\ \ \ 
   \mbox{for}\ \  j=1,2,\dots.
\end{equation}
Note that the eigenvalues are well separated. In particular, 
$$
   \frac{\lambda_{2,r_k}^{\ee}}{\lambda_{1,r_k}^{\ee}}=\frac1{3^{2r_k+2}}.
$$

For the Wiener case, $\lambda_{j,r_k}^{\ww}$'s are 
the eigenvalues of the operator

$$
(C_{1,r_k}^{\,\ww}f)(x)=\int_0^1K^{\ww}_{1,r_k}(x,t)\,f(t)\,\,{\rm  d}t.
$$
The eigenvalues $\lambda_{j,r_k}^{\ww}$ also satisfy the
Sturm-Liouville problem~\eqref{SL} but with different boundary
conditions
$$
f(0)=f^\prime(0)=\cdots=f^{(r_k)}(0)=f^{(r_k+1)}(1)=f^{(r_k+2)}(1)=
\cdots=f^{(2r_k+1)}(1)=0.
$$
The eigenvalues $\lambda_{j,r_k}^{\ww}$ are \emph{not}
exactly known. It is known \cite{GHT03} 
that they have the same asymptotic behavior
as in (\ref{Eulereig}),
\begin{equation}\label{Wiener-eig}
  \lambda_{j,r_k}^{\ww}=\left(\frac{1}{\pi(j-1/2)}\right)^{2r_k+2}
   +\OO\left(j^{-(2r_k+3)}\right)
   \ \ \ \ \  \mbox{as}\ \  j\to\infty.
\end{equation}
For tractability studies the asymptotic behavior is not enough
and the two largest eigenvalues play an essential role. 
That is why we will prove that
\begin{eqnarray*}
   \lambda_{1,r_k}^{\ww}&=&\frac1{(r_k!)^2}
   \left(\frac1{(2r_k+2)(2r_k+1)}+\OO(r_k^{-4})\right),\\
   \lambda_{2,r_k}^{\ww}&=&\Theta\left(\frac1{(r_k!)^2\,r_k^4}\right),
\end{eqnarray*}
where the factors in the big $\OO$ and $\Theta$ notations do not 
depend on $r_k$. 

Note that the largest eigenvalues for the Euler case go to zero
exponentially fast with~$r_k$, whereas for the Wiener case they go
to zero super exponentially fast due to the presence of factorials. 
However, the ratio of the two largest eigenvalues for the Wiener case,
$$
   \frac{\lambda_{2,r_k}^{\ww}}{\lambda_{1,r_k}^{\ww}}=\Theta(r_k^{-2}),
$$
is much larger than that for the Euler case.

\subsection{Tractability}

We present the precise definitions of four notions of tractability.
Let $n(\eps,d)$ denote the average case information complexity defined
in~\eqref{infocomp}, and let ${\rm APP}=\{{\rm APP}_d\}_{d=1,2,\dots}$ 
denote the sequence of multivariate approximation problems.  We say that
\begin{itemize}
\item ${\rm APP}$ is \emph{weakly tractable} iff
$$
  \lim_{\eps^{-1}+ d\to\infty}\frac{\ln\,n(\eps,d)}{\eps^{-1}+d}=0,
$$
with the convention that $\ln\,0=0$.
\item ${\rm APP}$ is \emph{quasi-polynomially tractable} iff there are
positive numbers $C$ and $t$ such that
$$
  n(\eps,d)\le C\,\exp\big(\,t\,(1+\ln\,d)\,(1+\ln\,\eps^{-1}\,\big)
  \ \ \ \ \ \mbox{for all}\ \ \ d=1,2,\dots, \ \ \eps\in(0,1).
$$
\item ${\rm APP}$ is \emph{polynomially tractable} iff there are 
non-negative numbers $C,q$ and $p$ such that
$$
  n(\eps,d)\le C\,d^{\,q}\,\eps^{-p}
  \ \ \ \ \ \mbox{for all}\ \ \ d=1,2,\dots, \ \ \eps\in(0,1).
$$ 
\item ${\rm APP}$ is \emph{strongly polynomially tractable} iff there 
are positive numbers $C$ and $p$ such that
$$
  n(\eps,d)\le C\,\eps^{-p}
  \ \ \ \ \ \mbox{for all}\ \ \ d=1,2,\dots, \ \ \eps\in(0,1).
$$ 
The infimum of $p$ satisfying the last bound is called the exponent
of strong polynomial tractability and denoted by $p^{\rm \,str-avg}$.
For the Euler and Wiener case, we use the notation $p^{\rm\,str-avg-x}$ with
$\xx\in\{{\rm E},{\rm W}\}$.
 \end{itemize}

Tractability can be fully characterized in terms of the eigenvalues 
$\lambda_{d,j}$. Necessary and sufficient conditions on weak, 
quasi-polynomial, polynomial and strong polynomial tractability 
can be found in Chapter~6 of~\cite{NW08} and Chapter~24 of~\cite{NW12} 
as well as in~\cite{LPW2} for non-homogeneous tensor products. 
For the Euler and Wiener integrated processes we need such conditions 
that are based on the sums of some power of the eigenvalues 
$\lambda_{d,j}$. We will cite these conditions when they are needed 
for specific tractability results.

\section{Euler Integrated Process}

We now analyze the Euler integrated process for which the eigenvalues
in the univariate cases are given by~\eqref{Eulereig}. Our aim is to
express tractability conditions in terms of the smoothness parameters
$\{r_k\}$.

\begin{thm}\label{Euler-thm}
Consider the approximation problem ${\rm APP}$ for the Euler
integrated process. 
\begin{itemize}

\item
${\rm APP}$ is weakly tractable iff 
\begin{equation}\label{weakEuler}
  \lim_{k\to\infty}r_k=\infty.
\end{equation}
Furthermore, if~\eqref{weakEuler} does not hold then we have the curse
of dimensionality since $n^{\ee}(\eps,d)$ depends exponentially on $d$ 
for each $\eps<1$.

\item ${\rm APP}$ is quasi-polynomially tractable iff 
\begin{equation}\label{Eulerquasi}
  \sup_{d\in\N}\ \frac1{\ln_+d}\ \sum_{k=1}^d(1+r_k)\,3^{-2r_k}<\infty,
\end{equation}
where $\ln_+d=\max(1,\ln\,d)$. 

\item ${\rm APP}$ is polynomially tractable iff ${\rm APP}$ is 
strongly polynomially tractable iff 
$$
  \sum_{k=1}^\infty 3^{-2\,\tau\,r_k}<\infty \ \ \ \ 
  \mbox{for some}\ \ \tau\in(0,1)
$$
or equivalently
iff
$$
   a_{\ee}:=\liminf_{k\to \infty}\ \frac{r_k}{\ln\,k}>\frac1{2\ln\,3}. 
$$
If so, then the exponent\footnote{It may happen that $a_{\ee}=\infty$.
Then the second term in the maximum defining $p^{\,\rm str-avg-E}$
is zero.} of strong polynomial tractability is
$$
  p^{\rm str-avg-E}
  =\max\left(\frac2{2r_1+1},\frac2{2a_{\ee}\,\ln\,3-1}\right).
$$
\end{itemize}

\end{thm}

We briefly comment on Theorem~\ref{Euler-thm}. First of all, we stress 
that polynomial and strong tractability are equivalent.
That is, these two notions coincide for the Euler integrated
process: in this case a \lq\lq weaker\rq\rq" property of 
polynomial tractability 
implies a \lq\lq stronger\rq\rq\ property of strong polynomial tractability.  
Weak tractability requires that the smoothness parameters 
$r_k$ go to infinity, however, the speed of
convergence is irrelevant. To obtain at least quasi-polynomial tractability, 
we need to assume that $r_k$ increases at least as $a_{\ee}\,\ln\,k$ 
with 
$a_{\ee} > 1/(2\ln\,3)$. Indeed, assume for simplicity that
$$
  a_{\ee}:=\lim_{k\to\infty}\frac{r_k}{\ln\,k}.
$$
exists. If $a_{\ee}<1/(2\,\ln\,3)$ then for any positive 
$\beta<1-2\,a_{\ee}\,\ln\,3$ we have
\begin{equation}\label{qpol111}
  n^{\ee}(\eps,d)\ge c_1(\beta)\,(1-\eps^2)\,
  \exp\left((c_2(\beta)\,d^{\,\beta}\right)
\end{equation}
for some positive functions $c_1$ and $c_2$ of $\beta$. Note 
that~\eqref{qpol111} contradicts quasi-polynomial tractability. 
The proof of~\eqref{qpol111} goes like follows. We will show later
that
$$
  n^{\ee}(\eps,d)
  \ge
  (1-\eps^2)\prod_{k=1}^d \left(1+3^{-2(r_k+1)}\right).
$$
Then each factor $1+3^{-2(r_k+1)}$ for large $j$  can be estimated from 
below by $\exp(-c(\beta)k^{-1+\beta})$. From this we easily 
obtain~\eqref{qpol111}.

If $a_{\ee}=1/(2\,\ln\,3)$ then we can have quasi-polynomial tractability 
as illustrated by an example of $\{r_k\}$ in the introduction.
Furthermore, for this example we do not have polynomial tractability.
However, it may also happen that for $a_{\ee}=1/(2\,\ln\,3)$ we do not have
quasi-polynomial tractability. For example, this is the case when
$$
   r_k=\left\lceil\frac{\ln_+k}{2\,\ln\,3}\right\rceil,
$$
which can be checked directly from~\eqref{Eulerquasi}.

On the other hand, if $a_{\ee}>1/(2\,\ln\,3)$ then we obtain strong
polynomial tractability. This shows that there is a ``thin'' zone 
of $\{r_k\}$ that separates quasi-polynomial and strong polynomial
tractabilities.

We now comment on the exponent of strong polynomial tractability. 
Note that for $a_{\ee}\ge (r_1+1)/\ln\,3$ we have
$$
   p^{\,\rm str-avg-E}=\frac2{2r_1+1}.
$$
In this case, the result is especially pleasing hence the complexity for 
any $d$ is roughly bounded by the complexity  for the univariate case. 
Furthermore, this happens for all $r_k$'s that tend to infinity faster 
than $\ln\,k$. On the other hand, if
$a_\ee\in(1/(2\,\ln\,3)),2(r_1+1)/(2\,\ln\,3))$ then we have 
$$
  p^{\,\rm str-avg-E}=\frac2{2a_{\ee}\,\ln\,3-1 },
$$
and $p^{\,\rm str-avg-E}$ can be arbitrarily large when $a_{\ee}$
is close to $1/(2\,\ln\,3)$.

\section{Wiener Integrated Process}

We now turn to the Wiener integrated process for which the eigenvalues
for the univariate cases $\lambda_{j,r_k}^{\ww}$ are only known 
asymptotically, see~\eqref{Wiener-eig}. To express tractability 
conditions in terms of the smoothness parameters $\{r_k\}$ we will need 
to prove the behavior of the two largest eigenvalues for large $r_k$.

\begin{thm}\label{Wiener-largest}
Consider the univariate Wiener process with the smoothness 
parameter $r$, and let $\lambda_{j,r}^{\ww}$'s 
denote the eigenvalues of 
the covariance operator $C_{1,r}^{\ww}$. Then

\begin{eqnarray*}
  \lambda_{1,r}^{\ww}&=&\frac1{(r!)^2}\left(\frac1{(2r+2)\,(2r+1)}
    +\OO(r^{-4})\right),\\
  \lambda_{2,r}^{\ww}&=&\Theta\left(\frac1{(r!)^2\,r^4}\right),\\
  \sup_{\tau\in[\tau_0,1]}\ 
  \frac{\sum_{j=3}^\infty\left[
  \lambda_{j,r}^{\ww}\right]^\tau}{\left[\lambda_{2,r}^{\ww}\right]^\tau}
&=&\OO(r^{-h})
  \ \ \ \ \ \ \ \ \   \mbox{for some}\ h>0\ \ \   
\mbox{and for all}\ \ \ \tau_0\in(\tfrac35,1].
\end{eqnarray*}
\end{thm}

Observe that the two largest eigenvalues for the Wiener case are much 
smaller than for the Euler case. On the other hand, their ratio for the 
Wiener case is much larger than for the Euler case. Therefore, 
the sequences $\{\lambda_{j,r}^{\ww}\}$ and $\{\lambda_{j,r}^{\ee}\}$
are quite different although they have the same asymptotic behavior. 

The uniform convergence in the last assertion of Theorem~\ref{Wiener-largest} 
at the neighborhood of $\tau=1$ is needed when we deal with 
quasi-polynomial tractability. The convergence for a specific $\tau$ is 
needed for strong polynomial and polynomial tractability. The lower bound 
$\tfrac35$ for $\tau_0$ is surely not sharp. A possible
improvement of this lower bound would improve the exponent of strong
polynomial tractability.

Based on the estimates presented in Theorem~\ref{Wiener-largest}
we will be able to express tractability conditions for the Wiener case 
in terms of $\{r_k\}$.

\begin{thm}\label{Wiener-thm}
Consider the approximation problem ${\rm APP}$ for the Wiener integrated 
process. 
\begin{itemize}

\item
${\rm APP}$ is weakly tractable iff 
\begin{equation}\label{weakWiener}
  \lim_{k\to\infty}r_k=\infty.
\end{equation}
Furthermore, if~\eqref{weakWiener} does not hold then we have the curse
of dimensionality since $n^{\ww}(\eps,d)$ depends exponentially on $d$ for
each $\eps<1$.

\item ${\rm APP}$ is quasi-polynomially tractable iff 
\begin{equation}\label{Wienerquasi}
  \sup_{d\in\N}\ \frac1{\ln_+d}\ \sum_{k=1}^d
  (1+r_k)^{-2}\,\ln_+r_k<\infty,
\end{equation}
where, we use $\ln_+x=\max(1,\ln\,x)$ for $x>0$, and $\ln_+0=1$. 
\item ${\rm APP}$ is polynomially tractable iff ${\rm APP}$ is strongly 
polynomially tractable iff 
$$
  \liminf_{k\to \infty}\ \frac{r_k}{k^s}>0
  \ \ \ \ \mbox{for some}\ \ \ \ \ s>\tfrac12.
$$
\end{itemize}
\end{thm}

We briefly comment on Theorem~\ref{Wiener-thm}. As for the Euler case, strong 
polynomial and polynomial tractability are equivalent, and weak tractability 
holds under the same condition $\lim_k r_k=\infty$. That ends the
similarity between the Wiener and Euler cases since the conditions on 
quasi-polynomial and polynomial tractability are quite different. For the 
Wiener case, we must assume that $r_k$'s go to infinity at least as fast 
as $k^{-s}$ for some $s>\tfrac12$. However, the zone between quasi-polynomial 
and polynomial tractabilities is again thin, as for the Euler case. 

It is worth to add that quasi-polynomial tractability plays a much
more important role in the worst case setting. 
The difference with the average case setting is due to the fact
that even for the constant sequence $r_k={\rm const}>0$ we have 
quasi-polynomial tractability in the worst case
setting as shown in~\cite{GW11}. 

We now discuss the exponent of strong tractability which is not addressed in 
Theorem~\ref{Wiener-thm}. For simplicity, let us assume that for some 
$s>\tfrac12$ we have 
$$
  r_k=k^s \ \ \ \ \mbox{for all}\ \ \ k\in\N.
$$
Then we have strong polynomial tractability and 
the exponent $p^{\rm str-avg-\ww}$ 
is given in~\eqref{expostpol} as the infimum of $2\tau/(1-\tau)$
for $\tau$ from $(0,1)$ which satisfies condition~\eqref{poltract} 
below with $q=0$.
From the proof of Theorem~\ref{Wiener-thm} we know that $\tau>1/(2r_1+2)$. 
Furthermore, \eqref{adstrong} implies that $\tau>1/(2s)$. 
These two estimates yield 
lower bounds on the exponent. On the other hand, our proof of strong 
polynomial tractability is valid only for $\tau>\frac35$, and this effects an 
upper bound on the exponent. Hence,
$$
  \max\left(\frac2{2r_1+1},\frac2{2s-1}\right)
  \le p^{\rm str-avg-\ww}
  \le \max\left(\frac2{2s-1},3\right).
$$
We stress that only for $s\in(\tfrac12,\tfrac56]$ we know the exponent exactly,
$
  p^{\rm str-avg-\ww}=\tfrac{2}{2s-1}.
$
Note that $p^{\rm str-avg-\ww}$ can be arbitrarily large if $s$ is close 
to $\tfrac12$.

For $s>\tfrac56$, our bounds on the eigenvalues $\lambda_{j,r_k}^{\ww}$ 
are too weak to get the exact value of the exponent but sufficient to deduce 
strong polynomial tractability.

\section{Proof of Theorem~\ref{Euler-thm}}

It is convenient to deal first with polynomial tractability. 
Let PT stand for polynomial tractability and SPT for strong polynomial
tractability. To prove this point of Theorem~\ref{Euler-thm} 
it is enough to show that
\begin{equation} \label{4stat}
a_\ee>\frac1{2\,\ln\,3} \ 
\Rightarrow\  
\sum_{k=1}^\infty 3^{-2\tau\,r_k}<\infty\ 
\Rightarrow \ 
{\rm SPT}\ 
\Rightarrow\ 
{\rm PT}
\Rightarrow \ 
a_\ee>\frac1{2\,\ln\,3} \ .
\end{equation}
The first claim,
$a_\ee>1/(2\,\ln\,3)\ \Rightarrow
S_\tau:=\sum_{k=1}^\infty3^{-2\tau\,r_k}<\infty$ for some
$\tau\in(0,1)$, is an easy calculus exercise.
Indeed, let $a_\ee>1/(2\ln\,3)$. Then for some $\delta>0$ and all $k$ large
enough we have
$\tfrac{r_k}{\ln k}> \tfrac{1+\delta}{2\ln 3}$, hence 
$3^{-2\tau r_k} < k^{-(1+\delta)\tau}$ and $S_\tau<\infty$ whenever
$\tfrac{1}{1+\delta}<\tau<1$. 

Recall now the polynomial tractability criteria.
We know from Chapter 6 of~\cite{NW08} that ${\rm APP}$ is polynomially 
tractable iff there exist $q\ge0$ and $\tau\in(0,1)$ such that
\begin{equation}\label{poltract}
   C:= \sup_{d\in\N}
   \frac{\left(\sum_{j=1}^\infty\lambda_{d,j}^\tau\right)^{1/\tau}}
   {\sum_{j=1}^\infty\lambda_{d,j}}\ d^{\, -q} <\infty.
\end{equation}
If so then
$$
  n(\eps,d)
  \le \left(\left(\frac{\tau\,C}{1-\tau}\right)^{\tau/(1-\tau)}+1\right)\,
      d^{\,q\,\tau/(1-\tau)}\,\eps^{-2\tau/(1-\tau)}
$$
for all $d\in\N$ and $\eps\in(0,1)$. 

Furthermore, ${\rm APP}$ is strongly polynomially tractable
iff~\eqref{poltract} holds with $q=0$. The exponent of strong
polynomial tractability is
\begin{equation}\label{expostpol}
   p^{\,\rm str-avg}=\inf\left\{\frac{2\tau}{1-\tau}\ \bigg| \ \ \tau \
  \mbox{satisfies~\eqref{poltract} with $q=0$}\right\}.   
\end{equation}

Motivated by condition (\ref{poltract}) and based on the explicit knowledge 
of the univariate eigenvalues for the Euler integrated 
process~\eqref{Eulereig}, 
we take $\tau\in(0,1)$ and obtain 
\begin{eqnarray*}
  \frac{\left(\sum_{j=1}^\infty\lambda_{d,j}^\tau\right)^{1/\tau}}
  {\sum_{j=1}^\infty\lambda_{d,j}}&=&
  \prod_{k=1}^d\frac{\left(\sum_{j=1}^\infty\left(
  \lambda_{j,r_k}^{\ee}\right)^\tau\right)^{1/\tau}}{\sum_{j=1}^\infty
  \lambda_{j,r_k}^{\ee}}\\
  &=&
  \prod_{k=1}^d\frac{\left(\sum_{j=1}^\infty(2j-1)^{-(2r_k+2)\tau}
  \right)^{1/\tau}}{\sum_{j=1}^\infty(2j-1)^{-(2r_k+2)}}\\
  &=&\prod_{k=1}^d
  \frac{\left(1+\sum_{j=2}^\infty(2j-1)^{-2\tau(r_k+1)}\right)^{1/\tau}}
  {1+\sum_{j=2}^\infty(2j-1)^{-2(r_k+1)}}.
\end{eqnarray*}
Since $r_k\ge r_1$, note that the expression above is finite for all 
$\tau\in(1/(2r_1+2),1)$. Furthermore for such $\tau$ we have 
$$
  3^{-2\tau(r_k+1)}
  \le\sum_{j=2}^\infty(2j-1)^{-2\tau(r_k+1)}
  \le 3^{-2\tau(r_k+1)}+\sum_{j=5}^\infty j^{-2\tau(r_k+1)},
$$
and
$$
  \sum_{j=5}^\infty j^{-2\tau(r_k+1)}
  \le \int_4^\infty x^{-2\tau(r_k+1)}\,{\rm d}x
  = \frac{4^{1-2\tau(r_k+1)}}{2\tau(r_k+1)-1}
  \le \frac{3}{2\tau(r_1+1)-1} \,3^{-2\tau(r_k+1)}.
$$
Therefore
\begin{equation}\label{useful}
  \frac{\left(\sum_{j=1}^\infty\lambda_{d,j}^\tau\right)^{1/\tau}}
  {\sum_{j=1}^\infty\lambda_{d,j}}=
  \prod_{k=1}^d\frac{\left(1 +a_k 3^{-2\tau(r_k+1)}\right)^{1/\tau}}
  {1+b_k 3^{-2(r_k+1)}},
\end{equation}
where $a_k\ge b_k$ and they are uniformly bounded, 
\begin{equation}\label{useful2}
  1\le a_k\le\frac{2\tau(r_1+1)+2}{2\tau(r_1+1)-1}\ \ \ \ \
  \mbox{and}\ \ \ \ \ 1\le b_k\le \frac{2r_1+4}{2r_1+1}.
\end{equation}

Assume now that $S_\tau<\infty$ for some $\tau<1$. By using
(\ref{useful}) and (\ref{useful2}) we obtain
\begin{eqnarray*}
\sup_{d}
\frac{\left(\sum_{j=1}^\infty\lambda_{d,j}^\tau\right)^{1/\tau}}
  {\sum_{j=1}^\infty\lambda_{d,j}}
&\le&   
 \prod_{k=1}^\infty \left(1 +a_k 3^{-2\tau(r_k+1)}\right)^{1/\tau}
\\
&\le&
\exp \left( \tau^{-1} \sup_k a_k \sum_{k=1}^\infty 3^{-2\tau(r_k+1)}\right)
\le
\exp \left( \tau^{-1} \sup_k a_k  \ S_\tau \right) <\infty.
\end{eqnarray*}
Hence, the criterion (\ref{poltract}) is verified 
with $q=0$, and we conclude that
$S_\tau<\infty$ $\Rightarrow$ {\rm SPT}.

Implication {\rm SPT} $\Rightarrow$ {\rm PT} is trivial.

Assume now that {\rm PT} holds. By (\ref{poltract}) and
(\ref{useful}) this implies that  
$$
  \prod_{k=1}^d
  \frac{\left(1 +a_k 3^{-2\tau(r_k+1)}\right)^{1/\tau}}
  {1+b_k 3^{-2(r_k+1)}}<C\,d^{\,q}
$$
for some $C,q\ge 0$ and $\tau\in (0,1)$. Moreover, is easy to check that
$$
  \frac{\left(1 +a_k 3^{-2\tau(r_k+1)}\right)^{1/\tau}}
  {1+b_k 3^{-2(r_k+1)}}\ge 1+c_k 3^{-2\tau(r_k+1)}
$$
for $c_k\ge a_k (1-3^{-2(r_k+1)(1-\tau)})/(1+b_k 3^{-2(r_k+1)})=\Omega(1)$.
Taking logarithms we conclude that 
$$
M:= \sup_d\ \frac1{\ln_+d}\ \sum_{k=1}^d  3^{-2\tau(r_k+1)}
<\infty.
$$
The sum with respect to $k$ can be lower bounded by $d \cdot 3^{-2\tau(r_d+1)}$,
as done at the beginning of the proof, and we obtain
$d \cdot 3^{-2\tau(r_d+1)}\le M \ln_+d$, which is equivalent to
\[
\frac{r_d+1}{\ln d} \ge 
\frac{1-\tfrac{\ln\ln_+d-\ln M}{\ln d}}{\tau\cdot2\ln 3},
\]   
and implies that $a_\ee\ge 1/(2\tau\,\ln\,3)  >1/(2\,\ln\,3)$, as claimed.
The equivalence of all statements in (\ref{4stat}) 
is therefore verified.

We now consider the exponent  $p^{\rm str-avg-E}$.
Assume now that $a_\ee>\tfrac{1}{2\,\ln\,3}$. 
Then, as already shown,  $\sum_{k=1}^\infty3^{-2\tau(r_k+1)}<\infty$
for all $\tau>\tfrac{1}{2a_\ee\ln\,3}$ and~\eqref{poltract} holds with $q=0$ if
$\tau>\tfrac{1}{2r_1+2}$. Hence, we obtain strong polynomial tractability.
Furthermore, $\tau$ can be taken in the limit as
$\tau_*:= \max\left(\tfrac{1}{2r_1+2},\tfrac{1}{2a_\ee\ln\,3}\right)$,
and~\eqref{expostpol} yields that the exponent of strong polynomial 
tractability is at most 
\[
  p_*:= \frac{2\tau_*}{1-\tau_*}
  = \max\left( \frac{2}{2r_1+1},\frac{2}{2a_\ee\ln\,3 -1}\right).
\]  
Conversely, assume that strong polynomial tractability holds. Then 
$$
  \prod_{k=1}^\infty
  \frac{\left(1 +a_k 3^{-2\tau(r_k+1)}\right)^{1/\tau}}
  {1+b_k 3^{-2(r_k+1)}}<\infty
$$
for some $\tau\in(0,1)$. Clearly, we must take $\tau>1/(2r_1+2)$ and 
$\tau>1/(2a_\ee\ln\,3)$. This implies that 
the exponent is at least $p_*$.
This completes the part of the proof related to polynomial and strong 
polynomial tractability.

We now turn to weak tractability. We know from~\cite{LPW2} that
${\rm APP}$ is weakly tractable if there exists $\tau\in(0,1)$ such
that 
\begin{equation}\label{weakcond}
  \lim_{d\to \infty}\
  \frac1d\ \sum_{k=1}^d\ \sum_{j=2}^\infty
  \left(\frac{\lambda_{j,r_k}^{\ee}}{\lambda_{1,r_k}^{\ee}}
   \right)^{\tau}=0.
\end{equation}
In our case, we have
$$
  \frac{\lambda_{j,r_k}^{\ee}}{\lambda_{1,r_k}^{\ee}}= (2j-1)^{-2(r_k+1)}.
$$
As before, for $\tau\in(\tfrac12,1)$ we have 
$$
  \sum_{j=2}^\infty(2j-1)^{-2\tau(r_k +1)}\le  
  \frac{2\tau(r_k+1)+2}{2\tau(r_k+1)-1}\,3^{-2\tau(r_k+1)}
  \le \frac{2(1+\tau)}{2\tau-1}\,3^{-2\tau(r_k+1)}.
$$
Assume that $\lim_{k\to\infty}r_k=\infty$. Then for an arbitrarily 
large $M$ there is an integer $k_M$ such that $r_k\ge M$ for all 
$k\ge k_M$. Hence, for $d\ge k_M$ we have 
$$
  \frac1d\ \sum_{k=1}^d\ \sum_{j=2}^\infty
  \left(\frac{\lambda_{j,r_k}^{\ee}}{\lambda_{1,r_k}^{\ee}}
  \right)^{\tau}
  \le    \frac{2(1+\tau)}{2\tau-1} \left(\frac{k_M}d+3^{-2\tau(M+1)}\right),
$$
and we obtain (\ref{weakcond}) by letting first $d$, 
and then $M$ go to infinity.

On the other hand, if $r=\lim_{k\to\infty}r_k<\infty$ 
then there is an integer $k_0$ such that $r_k=r$ for all $k\ge k_0$,
and the limit in~\eqref{weakcond} is not zero.
In this case, we prove that $n=n^{\ee}(\eps,d)$ is an
exponential function of $d$ 
and therefore weak tractability does not hold.
Indeed, we have
$$
  \sum_{j=1}^\infty\lambda_{d,j}-n\lambda_{d,1}\le 
  \sum_{j=n+1}^\infty\lambda_{d,j}
  \le \eps^2\sum_{j=1}^\infty\lambda_{d,j}, 
$$
and therefore
\begin{eqnarray*}
   n&\ge& (1-\eps^2)\sum_{j=1}^\infty\frac{\lambda_{d,j}}{\lambda_{d,1}}
   =(1-\eps^2)\left(\prod_{k=1}^{k_0-1}\sum_{j=1}^\infty
   \frac{\lambda_{j,r_k}^{\ee}}
   {\lambda_{1,r_k}^{\ee}}\right)
   \left(1+\sum_{j=2}^\infty
   \frac{\lambda_{j,r}^{\ee}}{\lambda_{1,r}^{\ee}} \right)^{d-k_0+1}
   \\
   &\ge&(1-\eps^2)\left(1+\sum_{j=2}^\infty
   \frac{\lambda_{j,r}^{\ee}}{\lambda_{1,r}^{\ee}} \right)^{d-k_0+1}.
\end{eqnarray*}
This bound is an exponential function of $d$. It contradicts weak 
tractability and completes the part of the proof related to this notion.

We finally consider quasi-polynomial tractability. 
We know from~\cite{LPW2} that ${\rm APP}$ is quasi-polynomially
tractable iff there exists a positive $\delta$ such that
\begin{equation}\label{quasicondition}
  \sup_{d\in\N}\ \frac{\sum_{j=1}^\infty\lambda_{d,j}^{1-\delta/\ln_+d}}
  {\left(\sum_{j=1}^\infty\lambda_{d,j}\right)^{1-\delta/\ln_+d}}<\infty,
\end{equation}
where $\ln_+d=\max(1,\ln\,d)$. 

Sufficiency. We first prove that~\eqref{Eulerquasi} 
implies~\eqref{quasicondition} with $\delta=\tfrac12$. Let
$$
  \lambda(j,k)=(2j-1)^{-2(r_k+1)}.
$$
We have
$$
  \sup_{d\in\N}\ \frac{\sum_{j=1}^\infty\lambda_{d,j}^{1-\frac1{2\ln_+d}}}
  {\left(\sum_{j=1}^\infty\lambda_{d,j}\right)^{1-\frac1{2\ln_+d}}}=
   \sup_{d\in\N}\ \prod_{k=1}^d 
  \frac{ \sum_{j=1}^\infty \lambda(j,k)^{1-\frac{1}{2\ln_+d}}}
  {\left(  \sum_{j=1}^\infty\lambda(j,k)\right)^{1-\frac1{2\ln_+d}}} .
$$
We split the last product into two products
\[
  \Pi_1(d) := \prod_{k=1}^d 
  \left(\sum_{j=1}^\infty \lambda(j,k)\right)^{\frac{1}{2\ln_+d}}  
\]
and
\[
   \Pi_2(d):=  \prod_{k=1}^d 
            \frac{ \sum_{j=1}^\infty \lam(j,k)^{1-\frac{1}{2\ln_+d}} }
            { \sum_{j=1}^\infty \lam(j,k)  }.
\]
In what follows we use $C$ to denote a positive number which is 
independent of $d$ and~$\{r_k\}$, and whose value may change for  
successive estimates. For $\Pi_1(d)$ we simply have
\begin{eqnarray*}
  \Pi_1(d) &=& \prod_{k=1}^d 
  \left(  1+ \sum_{j=2}^\infty \lam(j,k)   \right)^{\frac{1}{2\ln_+d}}  
  \le \exp\left( \frac{1}{2\ln_+d}   \sum_{k=1}^d  
  \sum_{j=2}^\infty \lam(j,k)  \right) 
  \\
  &\le& \exp\left( \frac C{\ln_+d}   \sum_{k=1}^d  \lam(2,k)  \right) 
  = \exp\left( \frac C{\ln_+d}   \sum_{k=1}^d  3^{-2(r_k+1)}  \right). 
\end{eqnarray*}
Clearly, \eqref{Eulerquasi} implies that $\sup_{d\in\N} \Pi_1(d) <\infty$.

We now turn to the product $\Pi_2(d)$. We estimate each of its factors by 
\begin{eqnarray} \nonumber
  \frac{ \sum_{j=1}^\infty \lam(j,k)^{1-\frac{1}{2\ln_+d}} }
            { \sum_{k=1}^\infty \lam(j,k)     } 
  &\le&
  \frac{ 1+ \lam(2,k)^{1-\frac{1}{2\ln_+d}}+\sum_{j=3}^\infty 
  \lam(j,k)^{1-\frac{1}{2\ln_+d}} } { 1+  \lam(2,k)     } 
  \\   \label{suf0}
  &\le& \frac{ 1+ \lam(2,k)^{1-\frac{1}{2\ln_+d}}  }{ 1+  \lam(2,k)     }  
  + \sum_{j=3}^\infty \lam(j,k)^{1-\frac{1}{2\ln_+d}} .            
\end{eqnarray}

Note that if $|\ln\lam(2,k)|\le 3\ln_+d$, then
\begin{eqnarray*}
   \frac{ 1+ \lam(2,k)^{1-\frac{1}{2\ln_+d}} } { 1+  \lam(2,k)} 
   &=&  \frac{ 1+\lam(2,k) 
   \exp\left({\frac{-\ln\lam(2,k)}{2\ln_+d}}\right) } { 1+\lam(2,k)} 
   \\
   &\le& \frac{ 1+\lam(2,k)\left({1+ \frac{C|\ln\lam(2,k)|}{\ln_+d}}\right)}
   {1+  \lam(2,k)} 
   \\
   &\le&  1+  \frac{C\lam(2,k) |\ln\lam(2,k)|} {\ln_+d},
\end{eqnarray*}
while if $|\ln \lam(2,k)|\ge 3\ln_+d$, then
\[
   \frac{ 1+ \lam(2,k)^{1-\frac{1}{2\ln_+d}} } { 1+  \lam(2,k)  }
   \le  1+ \lam(2,k)^{1-\frac{1}{2\ln_+d}} 
   \le 1+ \lam(2,k)^{1/2} 
   \le  1+ d^{-3/2}.
\]
Thus, in any case
\begin{equation}\label{suf1}
\frac{ 1+ \lam(2,k)^{1-\frac{1}{2\ln_+d}} } { 1+  \lam(2,k)     } 
\le 1+ d^{-3/2} + \frac{C\lam(2,k) |\ln\lam(2,k)|} {\ln_+d}.
\end{equation}
Next, we have
\begin{equation} \label{suf1a}
   \sum_{j=3}^\infty \lam(j,k)^{1-\frac{1}{2\ln_+d}}
   \le  C \lam(3,k)^{1-\frac{1}{2\ln_+d}}
   =  C \lam(2,k)^{\frac{\ln 5}{\ln 3} (1-\frac{1}{2\ln_+d })}.
\end{equation}
We now show that~\eqref{Eulerquasi} implies that 
$\lam(2,k)=3^{-2(r_k+1)}\le C/k$.
First of all note that~\eqref{Eulerquasi} implies that 
$\lim_k r_k=\infty$, so that only finitely many initial $r_k$ may be
zero. Assume that $d$ is so large that $r_d\ge 1$ and $d\ge 3$.  
Since $(1+r_k)3^{-2r_k}$ is non-increasing, we have 
$$
  r_d 3^{-2r_d}
  \le \frac1d\,\sum_{k=1}^d(1+r_k)3^{-2r_k}\le\frac{C\,\ln\,d}{d},
$$
so that $3^{2r_d}\ge 3^{2r_d}/r_d\ge d/(C\,\ln\,d)$ and  
$$
  r_d\ge \frac{\ln\,d-\ln(C\,\ln\,d)}{2\,\ln\,3}\ge C_1\,\ln\,d.
$$
Hence,
$$
  \lam(2,d)=3^{-2(r_d+1)} \le \frac{r_d3^{-2r_d}}{r_d}\le
  \frac{C\,\ln\,d}{r_d\,d}\le \frac{C}{C_1\,d}\,
$$
as claimed. By enlarging the constant, we obtain the same inequality for 
{\it all} $d$. For $k\le d$, we then have by (\ref{suf1a})
\begin{equation}\label{suf2}
   \sum_{j=3}^\infty\lambda(j,k)^{1-\frac1{2\ln_+d}}
   \le C\, k^{-\frac{\ln(5)}{\ln(3)}(1-\frac1{2\ln_+d})}
   \le C\,k^{-\frac{\ln(5)}{\ln(3)}}
\end{equation}

Using $1+x\le \exp(x)$, from (\ref{suf0}), (\ref{suf1}), and (\ref{suf2}),
we obtain
\[
    \frac{ \sum_{j=1}^\infty \lam(j,k)^{1-\frac{1}{2\ln_+d}} }
            { \sum_{j=1}^\infty \lam(j,k)     } 
    \le \exp\left(  d^{-3/2} + \frac{C\lam(2,k) |\ln\lam(2,k)|} {\ln_+d} 
        + C  k^{-\frac{\ln 5}{\ln 3}}  \right).
\]
Then it follows that
\begin{eqnarray*}
   \Pi_2(d) &\le& \exp\left ( \sum_{k=1}^d \left(  d^{-3/2} + 
   \frac{C\lam(2,k) |\ln\lam(2,k)|} {\ln_+d} 
    + C k^{-\frac{\ln 5} {\ln 3} } \right) \right)
   \\
   &\le& \exp\left(\sum_{k=1}^d \left(  d^{-3/2} 
     + \frac{C 3^{-2r_k} (r_k+1)} {\ln_+d} 
     + C  k^{-\frac{\ln 5}{\ln 3} }   \right)\right),
\end{eqnarray*}
and~\eqref{Eulerquasi} implies that  $\sup_{d\in\N} \Pi_2(d) <\infty$. 
Therefore, 
 \[
    \sup_{d\in\N} \Pi_1(d)\,\Pi_2(d) \le  \sup_{d\in\N} \Pi_1(d) \,
    \ \sup_{d\in\N} \Pi_2(d) <\infty,
\] 
the required property~\eqref{quasicondition} is verified, so that 
the quasi-polynomial tractability is proved.
\medskip

Necessity. \ Assume now that quasi-polynomial tractability  holds. 
We prove in~\cite{LPW2} that quasi-polynomial tractability 
implies 
\begin{equation}\label{qpolnec}
    \sup_{d\in\N}\, \frac1{\ln_+d}\,\sum_{k=1}^d\sum_{j=1}^\infty
    \frac{\lambda(j,k)}{\Lambda(k)}\,\ln\left(\frac{\Lambda(k)}
    {\lambda(j,k)}\right)<\infty,
\end{equation}
where $\Lambda(k)=\sum_{j=1}^\infty\lambda(j,k)$. Clearly,
$\Lambda(k)/\lambda(j,k)>1$ so that all terms in the sums over $j$ are
positive. We simplify the last condition by omitting all terms for
$j\not=2$, and obtain
\begin{equation}\label{conse1}
   \sup_{d\ge\N} \ \frac{1}{\ln_+d} \  
   \sum_{k=1}^d  \frac{\lam(2,k)}{\Lam(k)} \, 
   \ln\left(\frac{\Lam(k)} {\lam(2,k)} \right) <\infty.
\end{equation}
Next, since $\Lam(k)>1$ we can also omit $\ln\, \Lambda(k)$ and obtain
\[
   \sup_{d\in\N} \ \frac{1}{\ln_+d} \ 
   \sum_{k=1}^d  \frac{\lam(2,k)}{\Lam(k)} 
   \, \ln\left(\frac{1} {\lam(2,k)} \right) <\infty.
\]
Furthermore, since $\{\Lam(k)\}$ is non-increasing, we have 
\[
   \sup_{d\in\N} \ \frac{1}{\ln_+d} \ 
   \sum_{k=1}^d  \lam(2,k) \, \ln\left(\frac{1} {\lam(2,k)} \right) 
   <\infty.
\]
This is equivalent to~\eqref{Eulerquasi}, and completes the proof. 
\qed

\section{Proof of Theorem~\ref{Wiener-largest}}

We represent the $r$-times integrated Wiener
process $W_r$ through a white-noise integral representation
\begin{equation}\label{wr}
   W_r(t):= \int_0^1 \frac{(t-u)_+^r}{r!}\,{\rm d}\,W(u),
\end{equation}
where the integration is carried over a standard Wiener process $W$
defined over $[0,1]$.
Clearly,
\begin{eqnarray} \nonumber
  \E \|W_r\|_2^2&=&\sum_{j=1}^\infty\lambda_{j,r}^{\ww}=
  \int_0^1K_{1,r}^{\ww}(t,t)\,{\rm d}t
  \\ \label{upperlargest}
  &=&\int_0^1\left(\int_0^t\frac{(t-u)^{2r}}{r!^2}\,
  {\rm d}u\right)\,{\rm d}t=\int_0^1\frac{t^{2r+1}}{(2r+1)r!^2}\,{\rm d}t
  = \frac 1{(2r+2)(2r+1)r!^2}.
\end{eqnarray}

We now supply a lower bound on the sum $\sum_{j=2}^\infty\lambda_{j,r}^{\ww}$.
To do this, we approximate $W_r$ by 
\[
   V_{r,1}(t):= t^r \, W_r(1) = \frac1{r!} 
   \int_0^1 t^r (1-u)^r\,{\rm d}W(u)\ \ \ \ \
   \mbox{for all}\ \ \ t\in[0,1].
\]
The process $V_{r,1}$ is of rank $1$ since  
$V_{r,1}(t):=\xi_1(\omega) \psi_1(t)$, where $\psi_1(t)=t^r/r!$
and $\xi_1(\omega)= \int_0^1  (1-u)^r  dW(u)$.
We now prove the following lemma. 

\begin{lemma}\label{lemma1}
For any $r>1$ we have
\begin{equation}\label{WV_rt}
   \E |W_r(t)-V_{r,1}(t)|^2 \le 
   \frac1{r!}\,\frac{3r^2}{(2r-2)^3} \, t^{2r-2}  (1-t)^2
   \ \ \ \ \ \mbox{for all}\ \ \ t\in[0,1],
\end{equation}
and
\begin{equation}\label{WV_r}
   \E ||W_r-V_{r,1}||_2^2 \le  \frac1{r!^2}\,  \frac{6r^2}{(2r-2)^6}.
\end{equation}
\end{lemma}

Before we prove the lemma, we stress that the order of the right hand 
side in (\ref{WV_r}) is smaller than that  of $\E||W_r||_2^2$. This means 
that $V_{r,1}$ incorporates the essential part of $W_r$ for large~$r$.

\noindent{\bf Proof of Lemma~\ref{lemma1}.} \
Let $\ed{0\le u\le t}$ be the characteristic function of $[0,t]$,
i.e., $\ed{0\le u\le t}=1$ for $u\in [0,t]$ and
$\ed{0\le u\le t}=0$ for $u\notin [0,t]$. 
We have
\begin{eqnarray*}
  \E |W_r(t)-V_{r,1}(t)|^2 
  &=& \frac1{r!^2} \int_0^1 
  \left[  t^r(1-u)^r  -(t-u)^r\ed{0\le u\le t} \right]^2\,{\rm d}u 
  \\
  &=& \frac{t^{2r}}{r!^2} \int_0^t   (1-u)^{2r}  
  \left[  1 -\left( \frac{t-u}{t(1-u)}\right)^r  \right]^2\,{\rm d}u 
  + \frac{t^{2r}}{r!^2} \int_t^1  (1-u)^{2r}\,{\rm d}u 
  \\
  &=& \frac{t^{2r}}{r!^2} \int_0^t   (1-u)^{2r}  
  \left[  1 -\left( 1- \frac{(1-t)u}{t(1-u)}\right)^r  \right]^2\,{\rm d}u 
    + \frac{t^{2r}}{r!^2}\int_t^1  (1-u)^{2r}\,{\rm d}u 
\\
 &:=& 
\frac{t^{2r}}{r!^2} [ I_1+ I_2].
\end{eqnarray*}

For $I_1$, we use an elementary bound
$0\le 1-(1-h)^r\le rh$ and get

\begin{eqnarray*}
   I_1  &\le& \int_0^t (1-u)^{2r}  r^2  \frac{(1-t)^2u^2}
   {t^2(1-u)^2}\,{\rm d}u 
   \\
   &=&  r^2   (1-t)^2 t^{-2}   \int_0^t (1-u)^{2r-2}   u^2\,{\rm d}u 
   \\
   &\le&  r^2   (1-t)^2 t^{-2}   \int_0^\infty  \exp(-(2r-2)u)
   u^2\,{\rm d}u 
   \\
   &=&  \frac{2r^2}{(2r-2)^3}\   (1-t)^2 t^{-2}  . 
\end{eqnarray*}
On the other hand,
\[
   I_2=   \int_0^{1-t}  v^{2r}\,{\rm d}v =\frac{(1-t)^{2r+1}}{2r+1}
    \le  \frac{r^2}{(2r-2)^3}\   (1-t)^2 t^{-2}. 
\]
By summing up we obtain
\[
   \E |W_r(t)-V_{r,1}(t)|^2 
   \le \frac1{r!^2}\,\frac{3r^2}{(2r-2)^3}  t^{2r-2}  (1-t)^2 , 
\]
as claimed in the first estimate of the lemma.
The second claim is obtained by a simple integration:
\begin{eqnarray*}
   \E ||W_r-V_{r,1}||_2^2 &=& \int_0^1 \E |W_r(t)-V_{r,1}(t)|^2\,{\rm d}t
   \\
   &\le&   \frac1{r!^2}\,\frac{3r^2}{(2r-2)^3}  
   \int_0^1 t^{2r-2}(1-t)^2\,{\rm d}t
   \\
   &=&   \frac1{r!^2}\,\frac{3r^2}{(2r-2)^3}  
   \int_0^1   (1-t)^{2r-2}  t^2\,{\rm d}t
   \\
   &\le&   \frac1{r!^2}\,\frac{3r^2}{(2r-2)^3}  
   \int_0^\infty  \exp(-(2r-2)t)  t^2 \,{\rm d}t
   \\
   &=&   \frac1{r!^2}\,\frac{6r^2}{(2r-2)^6}.
\end{eqnarray*}
as claimed. \qed

From Lemma~\ref{lemma1} we conclude that
$$
  \sum_{j=2}^\infty   \lam_{j,r}^{\ww}= 
  \inf_{V\ {\rm is \ rank \ one}} \E||W_r-V||_2^2
  \le  \E||W_r- V_r||_2^2 \le \frac{C}{r!^2\,r^{4}}.
$$

This fact and~\eqref{upperlargest} yield
\[
   \lam_{1,r}^{\ww}= \frac{1}{r!^2}
   \left(  \frac{1} {(2r+2)(2r+1)}+O(r^{-4})\right),
\]
as claimed in Theorem \ref{Wiener-largest}.
\medskip

We now proceed to estimates on the second largest eigenvalue
$\lambda_{2,r}^{\ww}$ for large $r$. Obviously, 
\begin{equation}\label{uppertwoeig} 
  \lambda_{2,r}^{\ww}\le 
  \sum_{j=2}^\infty   \lam_{j,r}^{\ww}=
  \OO\left(\frac{1}{r!^2\,r^4}\right). 
\end{equation}
We now show that the last bound is essentially sharp.
To do this we approximate $W_r$ by 
\[
   V_{r,2}(t):= \frac1{r!}
\int_0^1 \left[ t^r (1-u)^r - rt^{r-1}(1-t) u(1-u)^{r-1}\right]
\,{\rm d}W(u)\ \ \ \ \ \mbox{for all}\ \ \ t\in[0,1].
\]
The process $V_{r,2}$ is of rank $2$ since 
\[
  V_{r,2}(t) =\xi_1(\omega) \psi_1(t) -r \xi_2(\omega) \psi_2(t), 
\]
 where
\begin{eqnarray*}
\xi_1(\omega)= \int_0^1  (1-u)^r\,{\rm d}W(u)\ \ \ &\mbox{and}&\ \ \
  \psi_1(t)= \frac{t^r}{r!},\\ 
\xi_2(\omega)= \int_0^1 u\, (1-u)^{r-1}\,{\rm d}W(u)\  \ \ 
&\mbox{and}&\ \ \   
\psi_2(t)=\frac{t^{r-1}(1-t)}{r!}.
\end{eqnarray*}
Note that the term $\xi_1 \psi_1 $ coming from rank 1 
approximation is dominating in the rank 2 approximation, since
\[
   \E\,\xi_1^2 ||\psi_1||_2^2 
=\int_0^1(1-u)^{2r}\,{\rm d}u \cdot \frac1{r!^2} 
\cdot \int_0^1 t^{2r} dt
= \frac1{r!^2} \frac{1}{(2r+1)^2} \approx \frac1{r!^2}\, r^{-2},
\]
while for the correction term  $r \xi_2 \psi_2$ we have 
\[
   r^2 \E\, \xi_2^2 ||\psi_2||_2^2 
=r^2   \int_0^1u^2(1-u)^{2r-2}\,{\rm d}u \cdot \frac1{r!^2} 
\cdot \int_0^1 t^{2r-2} (1-t)^2dt
\approx \frac1{r!^2} r^{-4}.
\]
A careful analysis shows that the second eigenvalue of the covariance operator 
of $V_{2,r}$ is also of order $\frac1{r!^2} r^{-4}$. In other words,
there exists a positive $C$ independent of $r$ such that
\begin{equation}\label{k2}
\inf_{V \ {\rm is \, rank \, one}} \E||V_{r,2}-V||_2^2
\ge \frac{C}{r!^2\,r^{4}}.
\end{equation}

We now estimate how well $V_{r,2}$ approximates $W_r$.

\begin{lemma}\label{lemma2}
For any $r>2$ we have
\begin{equation}\label{WV_rt2}
  \E |W_r(t)-V_{r,2}(t)|^2 \le \frac1{r!^2}\,
  \frac{14 r^2 (r-1)^2}{(2r-4)^{5}} \, t^{2r-4}  (1-t)^4
  \ \ \ \ \ \mbox{for all}\ \ \ t\in[0,1],
\end{equation}
and
\begin{equation}\label{WV_r2}
  \E ||W_r-V_{r,2}||_2^2 \le  \frac1{r!^2}\,
  \frac{24\cdot 14\cdot r^2(r-1)^2}{(2r-4)^{10}}
  =\OO \left(\frac1{r!^2\,r^6}\right).
\end{equation}
\end{lemma}

The proofs of~\eqref{WV_rt2} and \eqref{WV_r2}
repeat (mostly, but not entirely) line by line those 

of Lemma~\ref{lemma1} but we provide them for the sake of completeness. 
These proofs also clearly indicate how higher order approximations can 
be handled. As in Lemma~\ref{lemma1} we again stress that the the order 
of the right hand side in~\eqref{WV_r2} is smaller than the rank 1 
approximation error computed in  (\ref{WV_r}). Therefore, rank 2 
approximation $V_{r,2}$ performs much better than rank 1 approximation 
$V_{r,1}$ for approximation of $W_r$ when $r$ is large.

\noindent{\bf Proof of Lemma~\ref{lemma2}.} \ 
Let $a:=\E |W_r(t)-V_{r,2}(t)|^2$. We have
\begin{eqnarray*}
  a&=& \frac1{r!^2} \int_0^1 \left[  t^r(1-u)^r -  
  r t^{r-1}(1-t) u(1-u)^{r-1}   -(t-u)^r\ed{0\le u\le t}\right]^2\,{\rm d}u 
  \\
  &=& \frac{t^{2r}}{r!^2}\int_0^t   (1-u)^{2r}  
  \left[  1 - \frac{ r(1-t) u}  {t(1-u)}   - 
  \left( \frac{t-u}{t(1-u)}\right)^r  \right]^2\,{\rm d}u 
  \\
  && \quad  + \frac{t^{2r}}{r!^2}
  \int_t^1  \left((1-u)^{r}-\frac{r(1-t)u}{t} (1-u)^{r-1}\right)^2\,{\rm d}u 
  \\
  &=& 
  \frac{t^{2r}}{r!^2}  
  \int_0^t   (1-u)^{2r}  \left[ 1- 
  \frac{ r(1-t) u}  {t(1-u)} -\left(  1 - \frac{(1-t)u}{t(1-u)}\right)^r  
  \right]^2 {\rm d}u 
  \\
  && \quad +  \frac{t^{2r}}{r!^2}   \int_t^1 
  \left((1-u)^{r} -   \frac{r(1-t)u}{t} (1-u)^{r-1}
  \right)^2\,{\rm d}u =: \frac{t^{2r}}{r!^2}\,[ I_1+ I_2].
\end{eqnarray*}

For $I_1$, we use an elementary bound
$0\ge  1-rh - (1-h)^r\ge -\frac{r(r-1)}{2}h^2$ and get
\begin{eqnarray*}
   I_1  &\le& \int_0^t  (1-u)^{2r}    \left(\frac{r(r-1)}{2} \cdot 
   \frac{(1-t)^2u^2} {t^2(1-u)^2}\right)^2\,{\rm d}u 
   \\
   &=&  \frac{r^2(r-1)^2}{4}    (1-t)^4 t^{-4}   
   \int_0^t (1-u)^{2r-4}   u^4\,{\rm d}u 
   \\
   &\le&   \frac{r^2(r-1)^2}{4} (1-t)^4 t^{-4}   
   \int_0^\infty  \exp(-(2r-4)u) u^4\,{\rm d}u 
   \\
   &=&   \frac{6r^2(r-1)^2}{(2r-4)^5}      (1-t)^4 t^{-4}  . 
\end{eqnarray*}
On the other hand, we can give the following, rather crude,
estimate for $I_2$. Note that  for $u>t$ and $r>1$ we have
\[
   \frac{r(1-t)u}{t} (1-u)^{r-1}  = r\cdot  
   \frac{(1-t)u}{t(1-u)} \cdot (1-u)^{r} \ge   (1-u)^{r}.
\]
Therefore,
\begin{eqnarray*}
   I_2 &\le&   \int_t^1  
  \left(  \frac{r(1-t)u}{t} (1-u)^{r-1}  \right)^2\,{\rm d}u 
\\
&\le&    \frac{r^2(1-t)^4}{t^4}     \int_t^1    u^2
  (1-u)^{2r-4}\,{\rm d}u 
\\
&\le&    \frac{r^2(1-t)^4}{t^4}     \int_0^\infty   u^2  
\exp(-(2r-4)u)\,{\rm d}u 
\\
 &=&    \frac{2r^2(1-t)^4}{(2r-4)^3t^4} \le 
     \frac{8r^2(r-1)^2}{(2r-4)^5}      (1-t)^4 t^{-4}  . 
\end{eqnarray*}
By summing up, we obtain
\[
   \E |W_r(t)-V_{r,2}(t)|^2 
   \le \frac1{r!^2}\,\frac{14 r^2(r-1)^2}{(2r-4)^5} (1-t)^4 t^{2r-4} , 
\]
as claimed in the first estimate of the lemma.
The second claim is obtained by a simple integration:
\begin{eqnarray*}
  \E\, ||W_r-V_{r,2} ||_2^2 &=& \int_0^1\,
  \E |W_r(t)-V_{r,2}(t)|^2\,{\rm d}t
  \\
  &\le&   \frac1{r!^2}\,
  \frac{14r^2(r-1)^2}{(2r-4)^5}  \int_0^1 (1-t)^4 t^{2r-4}\,{\rm d}t
  \\
  &=&   \frac1{r!^2}\,
  \frac{14r^2(r-1)^2}{(2r-4)^5} \int_0^1 (1-t)^{2r-4}  t^4\,{\rm d}t
  \\
  &\le&   \frac1{r!^2}\,
  \frac{14r^2(r-1)^2}{(2r-4)^5}     
  \int_0^\infty  \exp(-(2r-4)t)t^4\,{\rm d}t
  \\
  &=&   \frac1{r!^2}\,\frac{24 \cdot 14 r^2 (r-1)^2} {(2r-4)^{10}}, 
\end{eqnarray*}
as claimed. \qed

From Lemma~\ref{lemma2} we easily estimate $\lam_{2,r}^{\ww}$.
Let $\zeta \eta_1:=\zeta(\omega)\eta_1(t)$ be the first term of 
Karhunen-Lo\`eve expansion for $W_r$. Then
\begin{eqnarray*}
   \frac{C}{r!^2\,r^{4}} &\stackrel{\textrm{by}\ (\ref{k2})} {\le}& 
   \E\, ||V_{r,2}- \zeta\eta_1 ||_2^2
   \\
   &=& \E\, ||(V_{r,2}-W_r) + (W_r-\zeta\eta_1)  ||_2^2
   \\
   &\le& 2\,\E\, ||V_{r,2}-W_r||_2^2 + 2\,
   \E\,||W_r-\zeta \eta_1)||_2^2
   \\
   &=& 2\,\E\, ||V_{r,2}-W_r||_2^2 + 2\lam_{2,r}^{\ww} 
       +2\sum_{j=3}^{\infty}  \lam_{j,r}^{\ww} . 
\end{eqnarray*}
Since $V_{r,2}$ is a process of rank 2, we also have
\begin{equation}\label{sum3le}
   \sum_{j=3}^{\infty}  \lam_{j,r}^{\ww}
   =  \inf_{V\ {\rm of \ rank\ two}} \E\,||W_r-V||_2^2
   \le \E\,||W_r-V_{r,2}||_2^2.
\end{equation}
For future use, we combine this with~\eqref{WV_r2} and get
\begin{equation}\label{sum3ler6}
     \sum_{i=3}^{\infty}  \lam_{i,r}^{\ww} \le \frac{C_1}{r!^2\,r^{6}}. 
\end{equation}
Furthermore, (\ref{sum3le}) immediately yields
\[
   \frac{C}{r!^2\,r^{4}} \le 4\E ||V_{r,2}-W_r||_2^2 + 2\lam_{2,r}
   \stackrel{\textrm{by}\ (\ref{WV_r2})} {\le}
   \frac{C_1}{r!^2\,r^{6}} + 2\lam_{2,r}^{\ww}.
\]
This provides a lower bound for $\lam_{2,r}^{\ww}$ and together 
with~\eqref{uppertwoeig} proves that
\begin{equation}\label{secondeig}
   \lam_{2,r}^{\ww}=\Theta\left(\frac1{r!^2\,r^{4}}\right),
\end{equation}
as claimed.
\medskip

We are ready to prove the last assertion of 
Theorem~\ref{Wiener-largest}. To simplify notation, let 
$\lambda_{j,r}=\lambda_{j,r}^{\ww}$. We split the series 
$\sum_{j=3}^\infty\lambda_{j,r}$ into two pieces - a long 
but finite initial part and a tail. Let $M>2$ and 
$\tau\in[\tau_0,1]$ with $\tau_0\in(\tfrac35,1]$. Consider 
the initial part including $j=3,4,\dots,\lceil r^M\rceil$. Using 
H\"older's inequality we obtain 
\begin{eqnarray} \nonumber
   \sum_{j=3}^{\lceil r^M\rceil} \lam_{j,r}^\tau 
   &\le&
   \left(\sum_{j=3}^{\lceil r^M\rceil} \lam_{j,r}\right)^\tau 
   \left(\sum_{j=3}^{\lceil r^M\rceil} 1 \right)^{1-\tau}
   \\  \nonumber
   &\stackrel{\textrm{by}\ (\ref{sum3ler6})} {\le}&
   \left(\frac{C_1}{r!^2\,r^{6}}  \right)^{\tau}    r^{M(1-\tau)}
   = 
   r^{-2\tau}  \left(\frac{C_1}{r!^2\,r^{4}} \right)^{\tau} r^{M(1-\tau)}
   \\ \nonumber
   &\stackrel{\textrm{by}\ (\ref{secondeig})} {\le}&
   C \lam_{2,r}^{\tau}    r^{-2\tau + M(1-\tau)}\le
   C \lam_{2,r}^{\tau}    r^{-2\tau_0 + M(1-\tau_0)}.
\end{eqnarray}
Since $C$ can be taken independent of $\tau$, for some $h>0$ we have 
$$
  \sup_{\tau\in[\tau_0,1]}\ 
  \frac{\sum_{j=3}^{\lceil r^M\rceil} 
\lam_{j,r}^\tau}{\lambda_{2,r}^\tau}=\OO(r^{-h}),
  \quad \textrm{as}\ r\to\infty, 
$$
as long as 
\begin{equation}\label{Mle}
    M<\frac{2\tau_0}{1-\tau_0}.
\end{equation}

For the tail estimation of the eigenvalue series
$\sum_{j=\lceil r^M\rceil+1}^\infty\lambda_{j,r}$ we use
{\it approximation numbers} (or {\it linear widths}, in other 
terminology).

We need to recall the definition and few basic properties which we will 
use in the sequel. Let $A:B_1\to B_2$ be a bounded linear operator acting 
between two Banach spaces. The approximation number $a_n(A)$ for $n\ge 1$ 
is defined as 
\begin{equation}
  a_n(A):= \inf\left\{\, \|A-A_n\|\ \ 
  \big|\ \   \ A_n:B_1\to B_2\ \ \ \mbox{with}\ \ \ 
  {\rm rank}(A_n)<n \ \right\}.
\end{equation}
The following properties of $a_n(A)$ are well known, see~\cite{Pi87}.
\begin{itemize}

\item the sequence $\{a_n(A)\}_{n\in\N}$ is non-increasing, 

\item for the adjoint operator $A^*$ we have
\begin{equation} \label{anstar}
    a_n(A)=a_n(A^*),
\end{equation}

\item multiplicative property: \ for $A_1:B_1\to B_2$ and $A_2:B_2\to B_3$ 
we have 
\begin{equation}\label{anmult}
   a_{n+m-1} (A_2A_1) \le a_{n}(A_2)\,a_{m}(A_1)\ \ \ 
   \ \ \mbox{for all}\ \ \ n,m\in \N,
\end{equation}

\item if $A:H\to H$ is a self-adjoint compact operator acting for a Hilbert
  space $H$ with the non-increasing eigenvalues $\{\lam_n\}$ then
\begin{equation}\label{anH}
     a_n(A)=\lam_n.
\end{equation}
\end{itemize}

We will study approximation numbers for integration operators.
Let $I:L_2[0,1]\to L_2[0,1]$ be the conventional integration operator
\[
     (Ix)(t):=\int_0^t x(s)\,{\rm d}s\ \ \ \ \ \mbox{for all}
      \ \ \ t\in[0,1].
\]
Let $I^r$ denote the $r$-th iteration of $I$ for $r\ge1$. It is easy to 
check by induction that 
\begin{eqnarray*}
   (I^rx)(t)&=&\int_0^t\frac{(t-s)^{r-1}}{(r-1)!}\,x(s)\,{\rm d}s
    \ \ 
    \ \ \ \ \ \mbox{for all}\ \ \ t\in[0,1],\\
    ([I^r]^*x)(t)&=&\int_t^1\frac{(s-t)^{r-1}}{(r-1)!}\,
    x(s)\,{\rm d}s \ \ \ \ \ \ \ \mbox{for all}\ \ \ t\in[0,1],\\
    (I^r\,[I^r]^*)(t)&=&\int_0^1\left(
    \int_0^{\,\min(s,t)}
    \frac{\ (s-u)_+^{r-1}}{(r-1)!}\,\frac{\ (t-u)_+^{r-1}}{(r-1)!}\,{\rm d}u
    \right)\,x(s)\,{\rm d}s
    \ \ \ \ \ \mbox{for all}\ \ \ t\in[0,1].
\end{eqnarray*}
This shows that
$$
  C_{1,r}^{\ww}=I^{r+1}\,(I^{r+1})^*.
$$ 
We are interested in the approximation numbers of $I^r$. For $r=0$, it is 
well known that for some  positive $C$ we have 
\begin{equation}\label{anI1}
    a_n(I) \le C\, n^{-1}\ \ \ \ \mbox{for all}\ \ \ n\in\N,
\end{equation}
see \cite{ET}, pp. 118--119. 
We will extend this estimate for $I^r$ with an arbitrary $r$. Although
the constant we get is certainly not optimal, it suffices for our needs.

\begin{lemma}\label{lemma3}
We have 
\begin{equation}\label{anIr}
    a_n(I^r) \le C^{\,r}\, (2r)^{2r}\,n^{-r}\ \ \ \ \ 
    \mbox{for all}\ \ \ n,r\in\N,
\end{equation}
where $C$ is a constant from~\eqref{anI1}.
\end{lemma}

\noindent{\bf Proof of Lemma~\ref{lemma3}.}\ Let 
\[
   B_p := 2^{\,p \, 2^p}\ \ \ \ \ \mbox{for all}\ \ \ p=0,1,2,\dots\,.
\]
We will first prove by induction on $p$ that for any integer $p\ge 0$ 
we have 
\begin{equation} \label{anind}
     a_n(I^r) \le C^{\,r}\,B_p\, n^{-r} \ \ \ \ \ 
    \mbox{for all}\ \ \  n\ge 1\ \ \mbox{and}\ \ r\in[2^{p-1}, 2^{p}].
\end{equation}
For $p=0$ this fact is equivalent to~\eqref{anI1}. Assume that 
(\ref{anind}) holds for some integer~$p$. Take any integer 
$r\in [2^p,2^{p+1}]$ and write it as $r=r'+r''$ with
$2^{p-1}\le r_1,r_2\le 2^p$. By using $I^r=I^{r_1} I^{r_2}$ and the 
multiplicative property~\eqref{anmult}, we get for an odd index $2n-1$
\begin{eqnarray*}
  a_{2n-1} (I^r)  &=&  a_{2n-1} (I^{r_1}  I^{r_2} ) 
  \le  a_n(  I^{r_1})\,  a_n(I^{r_2})
  \\
  &\le&  C^{r_1} B_p n^{-r_1} \cdot  C^{r_2} B_p n^{-r_2}
  =  C^{r}\,B_p^2\, n^{-r} 
  \\
  &=& C^r\,B_p^2\,  2^r\, (2n)^{-r} \le  C^r\, 
  [ B_p^2  2^{2^{p+1} } ]\, (2n)^{-r} 
  \\
  &=&  C^r  2^ {\,2 p 2^p +2^{p+1}}\, (2n)^{-r} 
  =   C^r \, 2^ {(p+1) 2^{p+1}}\,  (2n)^{-r} 
  \\
  &=&    C^r \, B_{p+1} \ (2n)^{-r} \le C^r B_{p+1} \ (2n-1)^{-r}. 
\end{eqnarray*}
For an even index $2n$ we simply have 
\[
  a_{2n} (I^r) \le a_{2n-1} (I^r) \le  C^r  B_{p+1} \  (2n)^{-r} .
\]
Therefore, (\ref{anind}) is proved by induction.

For $r$ and $p$ as in (\ref{anind}), we have 
$B_p= (2^p)^{2^p}\le (2r)^{2r}$.
Hence, (\ref{anIr}) follows from (\ref{anind}).
\qed
\bigskip

We now relate approximation numbers $a_n(I^r)$ to the eigenvalues 
$\lambda_{j,r}$ of the operator $C_{1,r}^{\ww}=I^{r+1}(I^{r+1})^*$. 
We have 
\begin{eqnarray*}
  \lam_{2j,r} &\le& \lam_{2j-1,r}  
  \stackrel{\textrm{by}\ (\ref{anH})} {=}
  a_{2j-1} (I^{r+1} (I^{r+1})^{*})
  \\
  &\stackrel{\textrm{by}\ (\ref{anmult})} {\le}&
  a_j(I^{r+1}) a_j((I^{r+1})^{*})
  \stackrel{\textrm{by}\ (\ref{anstar})} {=}  a_j(I^{r+1})^2
  \\
  &\stackrel{\textrm{by}\ (\ref{anIr})} {\le}&
  C^{2(r+1)} (2(r+1))^{4(r+1)} j^{-2(r+1)}.
\end{eqnarray*}
This can be written as
\[
   \lam_{j,r}\le C_1^{\,r} \,r^{4(r+1)}\, j^{-2(r+1)}\ \ \ \ \ 
   \mbox{for all}\ \ \ r,j\in\N.
\]

Take a (small) positive $\alpha$. Consider
$r$ so large that $r\ge C_1^{1/\alpha}$ and $2(r+1)\tau>1$.
Then again for $\tau\in[\tau_0,1]$ we can sum up 
\begin{eqnarray*}
   \sum_{j=\lceil r^M\rceil +1}^\infty \lam_{j,r}^\tau  
   &\le& C_1^{\,r\tau}\,r^{4(r+1)\tau}  
   \sum_{j=\lceil r^M\rceil +1}^\infty  j^{-2(r+1)\tau}
   \\
   &\le&  r^{(4+\alpha)r\tau+4\tau}\int_{r^M}^\infty
   x^{-2(r+1)\tau}\,{\rm d}x\\
   &=& \frac{r^{r\tau(4+\alpha-2M)}\ r^{4\tau+M(1-2\tau)}}{2(r+1)\tau-1}.
\end{eqnarray*}
We relate the last estimate to $\lambda_{2,r}=\Theta(1/(r!^2r^4))$.
Since $r!=r^{r+1/2}\,e^{-r}\,\sqrt{2\pi}(1+o(1))$ by Stirling's formula, 
we have
$$
  \lambda_{2,r}=\frac{e^{2r}}{2\pi\,r^{2r+5}}\,(1+o(1)) \ \ \ \ \ \mbox{as}
  \ \ \ r\to\infty.
$$
Therefore
$$
  \sum_{j=\lceil r^M\rceil +1}^\infty\lambda_{j,r}^\tau=
  \OO\left(\lambda_{2,r}^\tau\,
  \frac{r^{-r\tau(2M-6-\alpha)}\ r^{9\tau+M(1-2\tau)}}
  {2(r+1)\tau-1}\,e^{-2r\tau}\right)
  =\OO\left(\lambda_{2,r}^\tau\,r^{-r\tau(2M-6-\alpha)}\right),
$$
where the factors in the big $\OO$ notation are independent
of $r,\tau$ and $\alpha$.

Assume that $M>3$. Then we can take a positive $\alpha$ such that 
$2M-6-\alpha>0$ and get
$$
  \sup_{\tau\in[\tau_0,1]}\  \frac{\sum_{j=\lceil r^M\rceil +1}^\infty
  \lambda_{j,r}^\tau}{\lambda_{2,r}^\tau}=\OO(r^{-h}),
  \quad \textrm{as}\ r\to\infty.
$$
Hence, 
$$
   \sup_{\tau\in[\tau_0,1]}\
  \frac{\sum_{j=1}^\infty \lambda_{j,r}^\tau}{\lambda_{2,r}^\tau}
  =\OO(r^{-h})  \quad \textrm{as}\ r\to\infty, 
$$
assuming \eqref{Mle} holds for some $M>3$. It is easy to
see that such a number $M$ exists since $\tau>\tfrac35$.
This completes the proof. \qed

\section{Proof of Theorem~\ref{Wiener-thm}}
As in the Euler case, we begin with polynomial tractability.
We now need to show that
$$
{\rm PT}\ \Rightarrow\ \liminf_k\frac{r_k}{k^s}>0\ \Rightarrow
\ {\rm SPT}\ \Rightarrow \ {\rm PT}.
$$

Observe that for $\lambda_{d,j}=\lambda_{d,j}^{\ww}$ and
$\tau\in(0,1)$, the expression in~\eqref{poltract} is now  
\begin{equation} \label{ad}
   a_d:=\frac{\left(\sum_{j=1}^\infty\lambda_{d,j}^\tau\right)^{1/\tau}}
        {\sum_{j=1}^\infty\lambda_{d,j}}
       = \prod_{k=1}^d
         \frac{\left(1+(\lambda_{2,r_k}/\lambda_{1,r_k})^\tau
   +\sum_{j=3}^\infty(\lambda_{j,r_k}/\lambda_{1,r_k})^\tau\right)^{1/\tau}}
         {1+\lambda_{2,r_k}/\lambda_{1,r_k}+\sum_{j=3}^\infty
           \lambda_{j,r_k}/\lambda_{1,r_k}}. 
\end{equation}
Since $\lambda_{j,r_k}=\Theta(j^{-2(r_k+1)})$ as $j\to \infty$,  
with the factors in the $\Theta$ notation depending on $r_k$,
then $a_d$ is finite 
iff $2(r_k+1)\tau>1$ for all $r_k$. Then $r_k\ge r_1$ implies 
that we need to consider $\tau\in( \tfrac{1}{2r_1+2},1)$.

Assume that we have polynomial tractability. Then 
$a_d\le C\,d^{\,q}$. 
Each ratio in the  product (\ref{ad}) is strictly larger than one. This
implies that $\lim_{k\to\infty}r_k=\infty$. 

Note that we can estimate $a_d$ from below by dropping the sums over $j$. 
Then
$$
  \prod_{k=1}^d
  \frac{\left(1+(\lambda_{2,r_k}/\lambda_{1,r_k})^\tau\right)^{1/\tau}}
  {1+2\lambda_{2,r_k}/\lambda_{1,r_k}}<C\,d^{\,q}\ .
$$
Taking logarithms and using 
the asymptotic formulas for $\lambda_{1,r_k}$ and $\lambda_{2,r_k}$ from
Theorem~\ref{Wiener-largest} yield
$$
 \sup_d\, \frac1{\ln_+d}\ \sum_{k=1}^d r_k^{-2\tau}<\infty.
$$
Since $d\,r_d^{-2\tau}\le \sum_{k=1}^dr_k^{-2\tau}$ we get 
$r_d^{-2\tau}=\OO(d^{-1}\,\ln_+d)$ and there exists $\delta>0$ such that 
$$
   r_d\ge \delta\,\left(\frac{d}{\ln_+d}\right)^{1/(2\tau)}\ \ \ \ \  
   \mbox{for all}\ \ \ d\in\N.
$$
Letting $s\in(\tfrac 12,\tfrac{1}{2\tau})$ we obtain 
\begin{equation}\label{rWiener}
   \liminf_{k\to\infty}\ \frac{r_k}{k^s}>0,
\end{equation}
as claimed. 

Assume now that~\eqref{rWiener} holds for some $s>\tfrac 12$.   
For $\tau\in(\max(\tfrac 35,\tfrac{1}{2s}),1]$ we can use the last assertion of
Theorem~\ref{Wiener-largest} to conclude that
\begin{equation}\label{adstrong}
   \sup_d a_d =\prod_{k=1}^\infty
  \frac{\left(1+\OO(r_k^{-2\tau})\right)^{1/\tau}} {1+\OO(r_k^{-2})}
  \le \exp\left\{ \OO\left(\sum_{k=1}^\infty r_k^{-2\tau}\right)\right\}
  = \exp\left\{\OO\left(\sum_{k=1}^\infty k^{-2s\tau}\right)\right\}<\infty.
\end{equation}
By criterion (\ref{poltract}) this implies strong polynomial and obviously 
polynomial tractability.

We turn to weak tractability.  Assume that
$\lim_{k\to\infty}r_k=\infty$. We verify the analogue of ~\eqref{weakcond}
for $\tau\in(\tfrac 35,1)$. From Theorem~\ref{Wiener-largest} we have 
$$
  b_d:=\frac1d\,\sum_{k=1}^d\sum_{j=2}^\infty\left(\frac{\lambda_{j,r_k}^{\ww}}
  {\lambda_{1,r_k}^{\ww}}\right)^{\tau}
  =\frac1d\,\sum_{k=1}^d\OO\left(r_k^{-2\tau}\right)
  =\OO \left(\frac1d\,\sum_{k=1}^d r_k^{-2\tau}\right).
$$
Clearly, $\lim_k r_k^{-2\tau}=0$ implies $\lim_d b_d=0$, which yields 
weak tractability. 

Let $r=\lim_{k\to\infty}r_k<\infty$. Then proceeding exactly as for the
Euler case, we can show that $n^{\ww}(\eps,d)$ is an exponential
function of $d$ which contradicts weak tractability and completes this
part of the proof. 

We finally consider quasi-polynomial tractability.
The proof is similar to the proof for the Euler case and
we only sketch it.  
We need to study~\eqref{quasicondition} and~\eqref{qpolnec}
for the Wiener eigenvalues.  
For~\eqref{quasicondition} we take $\delta=\tfrac12$ and 
$\tau_0\in(\tfrac35,1)$.
Let us chose $d_0$  such that $1-\tfrac{1}{2\ln\,d_0} \in[\tau_0,1]$. Then 
for all such $d\ge d_0$ we have $\tau_d:=1-1/(2\ln\,d)\in[\tau_0,1]$ 
and we can use the result on the uniform convergence presented in the 
last assertion of Theorem~\ref{Wiener-largest} with respect now to $d$.
Let denote
$Q_k:= \frac{\lambda_{2,r_k}}{\lambda_{1,r_k}}$.
We obtain
\begin{eqnarray*}
   \frac{\sum_{j=1}^\infty\lambda_{d,j}^{1-\delta/\ln_+d}}
   {\left(\sum_{j=1}^\infty\lambda_{d,j}\right)^{1-\delta/\ln_+d}}
   &=& \prod_{k=1}^d
   \frac{1+Q_k^{1-\frac1{2\,\ln\,d}}+
   \sum_{j=3}^\infty\left(\frac{\lambda_{j,r_k}}{\lambda_{1,r_k}}\right)^{1-
   \frac1{2\,\ln\,d}}}
   {\left(1+Q_k+
   \sum_{j=3}^\infty\frac{\lambda_{j,r_k}}{\lambda_{1,r_k}}\right)^{1-
   \frac1{2\ln\,d}}}
   \\
   &\le&\OO(1)\ \prod_{k=d_0}^d
   \frac{1+Q_k^{1-\frac1{2\,\ln\,d}} \left(1+o(r_k^{-h})\right)}
   {\left(1+Q_k\right)^{1-\frac1{2\ln\,d}}},
\end{eqnarray*}
with absolute constants as pre-factors in the $\OO(\cdot)$ notation.

Suppose that~\eqref{Wienerquasi} holds. 
Then $\lim_k r_k=\infty$ and 
$$
\prod_{k=d_0}^d \left(1+Q_k\right)^{\frac1{2\ln\,d}}
\le
\exp\left(\frac{2}{\ln\,d}\,\sum_{k=d_0}^d Q_k \right)
\le
\exp\left(\frac{C}{\ln\,d}\,\sum_{k=d_0}^d r_k^{-2}\right)
$$
is uniformly bounded in $d$. The factor 
$\prod_{k=d_0}^d
   \frac{1+Q_k^{1-\frac1{2\,\ln\,d}} \left(1+o(r_k^{-h})\right)}
   {1+Q_k}
$
can be analyzed exactly as for the Euler case.
By using $Q_k=\Theta (r_k^{-2})$, we have  
$$
   \frac{1+Q_k^{1-\frac1{2\,\ln\,d}} \left(1+o(r_k^{-h})\right)}
   {1+Q_k}
\le 
1+d^{-3/2}+C\,(1+r_k)^{-2} \left(\frac{\,\ln_+r_k}{\ln\,d}+o(r_k^{-h})\right).
$$
Recall that assumption~\eqref{Wienerquasi} 
yields $r_k^{-2}=\OO(\frac{\ln k}{k})$, hence
$$
\prod_{k=d_0}^d 
   \frac{1+Q_k^{1-\frac1{2\,\ln\,d}} \left(1+o(r_k^{-h})\right)}
   {1+Q_k}
\le \exp\left(\sum_{k=d_0}^d \left(d^{-3/2}+C\,(1+r_k)^{-2}
\left[ \frac{\ln_+r_k}{\ln\,d} +  r_k^{-h}\right]  \right)\right)
$$
is also uniformly bounded in $d$.
This means that~\eqref{Wienerquasi} implies quasi-polynomial tractability.

Suppose now that quasi-polynomial tractability holds.
Then we use~\eqref{qpolnec} and its consequence~\eqref{conse1},
which is equivalent to~\eqref{Wienerquasi}.
This completes the proof. \qed

\section*{Acknowledgment}

The work of the first and the third authors was done while they 
participated in the Trimester Program ``Analysis and 
Numerics for High Dimensional Problems'', May-August 2011, in Bonn, Germany, 
and enjoyed warm hospitality 
of the Hausdorff Research Institute for Mathematics. 

The work of the first author was supported by RFBR grants 
10-01-00154à, 11-01-12104-ofi-m, and by 
Federal Focused Program 2010-1.1.-111-128-033.
The work of the second and third authors was partially supported by
the National Science Foundation.

\medskip

\noindent{\bf Authors' Addresses:}

\vskip 1pc

\noindent
M. A. Lifshits, 
Department of Mathematics and Mechanics,\\
St. Petersburg State University,
198504 St. Petersburg, Russia,\\
email: mikhail@lifshits.org
\smallskip

\noindent
A. Papageorgiou,
Department of Computer Science,
Columbia University,\\
New York, NY 10027, USA, email: ap@cs.columbia.edu
\smallskip

\noindent
H. Wo\'zniakowski, Department of Computer Science, Columbia
University,\\
New York, NY 10027, USA, and \\
Institute of Applied Mathematics and Mechanics, University of Warsaw,\\
ul. Banacha 2, 02-097 Warszawa, Poland, email:
henryk@cs.columbia.edu
\vfill

\end{document}